\documentclass[10pt]{article}
\usepackage{amssymb}
\usepackage{bbm,booktabs,graphicx,multirow,epsfig}
\advance\hoffset by -1.8cm
\advance\voffset by -1.9cm 
\newcommand{\beq}{\begin{equation}}
\newcommand{\eeq}{\end{equation}}
\newcommand{\beqa}{\begin{eqnarray}}
\newcommand{\eeqa}{\end{eqnarray}}
\newcommand{\beqas}{\begin{eqnarray*}}
\newcommand{\eeqas}{\end{eqnarray*}}

\newcommand{\ceq}{\!\!\! & = & \!\!\!}
\newcommand{\cleq}{\!\!\! & \le & \!\!\!}

\newcommand{\real}{{\ensuremath{\mathbb{R}}}}

\newcommand{\el}{\epsilon}
\newcommand{\sumnl}{\sum\nolimits}

\newcommand{\limnl}{\lim\nolimits}
\newcommand{\infnl}{\inf\nolimits} 
\newcommand{\supnl}{\sup\nolimits}

\newcommand{\flx}{\lfloor x \rfloor}

\newcommand{\snx}{\sum\nolimits_{n\le x}}
\newcommand{\snn}{\sum\nolimits'_n}

\newcommand{\done}{\hfill $\Box$}

\newcommand{\I}{\mathcal{I}}

\newcommand{\R}{\mathrm{Res}}

\newtheorem{thm}{Theorem}[section]
\newtheorem{prop}{Proposition}[section]
\newtheorem{cor}{Corollary}[section]

\begin{document}

\begin{center}
\begin{Large}
{\bf On certain Gram matrices and their associated series}
\end{Large}
\vskip 6mm
\begin{large}
{\sl Werner Ehm}
\end{large}
\end{center}

\begin{abstract}
We derive formulae for Gram matrices arising in the Nyman--Beurling reformulation of the Riemann hypothesis. The development naturally leads upon series of the form $S(x) = \sum_{n\ge 1} R(nx)$ and their reciprocity relations. We give integral representations of these series; and we present decompositions of the quadratic forms associated with the Gram matrices along with a discussion of the components' properties. 
\end{abstract}

\section{Introduction}\label{intro}

In complex variable terms, the Nyman--Beurling criterion for the Riemann hypothesis [RH] (cf.~\cite{N,Be2}) amounts to the condition that the Mellin transform, $1/s$, of the indicator function of the unit interval can be approximated in $L^2(\Re(s) = 1/2, |ds|)$ by functions of the form $\sum_{n \le N} a_n \theta_n^s\,  \zeta(s)/s$, where $\zeta(s)$ is the Riemann zeta function, $\theta_n \in (0,1]$, and the coefficients $a_n$ may be taken to be real--valued. The criterion was significantly strengthened by B\'aez--Duarte \cite[Theorem 1.3]{BD1} who showed that one may take $\theta_n = 1/n$, so that it suffices to consider Dirichlet  polynomials times $\zeta(s)/s$ as the approximating functions. Let 
\beq
d_q^2(N) \equiv d_q^2(N;\, a_1,\ldots,a_N) = \frac{1}{2\pi}\int_{\Re(s) = 1/2} \bigg|\,1 -  \sumnl_{n\le N} a_n\, n^{-s}\,\zeta(s)\, \bigg|^2 \frac{|ds|}{|s|^{2q}}\quad (q = 1,2) \label{cplxerr}.
\eeq
For $q=1$, $d_q^2(N)$ represents the (squared) approximation error in B\'aez--Duarte's criterion which tells us that {\em RH is true if and only if}
\beq\label{BDcrc}
\liminf\nolimits_{N \to \infty}\,  \infnl_{a_1,\ldots,a_N}\, d_q^2(N;\, a_1,\ldots,a_N) = 0\, .
\eeq
The case $q=2$ will be dealt with in parallel as a variation of B\'aez--Duarte's criterion. Initially, the equivalence to RH is only known for $q=1$, yet as is easily seen, it likewise holds for $q=2$. 
A potential merit of the choice $q=2$ might lie in the stronger downweighting of the squared error in  (\ref{cplxerr}). 

On expanding the square and noting that $\bar s = 1-s$ on the crictical line  one obtains the decomposition 
\beq\label{qdec}
d_q^2(N;\, a_1,\ldots,a_N) = C^{(q)} -  2\sumnl_{n\le N} a_n F_n^{(q)} +  \sumnl_{m,n\le N} a_m a_n\,G_{m,n}^{(q)}
\eeq
where for $q = 1,2$
\beqa\label{Cdef}
C^{(q)} \ceq \frac{1}{2\pi i}\int_{\Re(s) = 1/2}\, \frac{ds}{s^q(1-s)^q} = q\, ,\\
F_n^{(q)} \ceq \frac{1}{2\pi i}\int_{\Re(s) = 1/2}\, n^{-s}\,\frac{\zeta(s)}{s^q(1-s)^q}\, ds\, , \label{Fdef}\\ 
G_{m,n}^{(q)} \ceq \frac{1}{2\pi i}\int_{\Re(s) = 1/2}\, m^{-s}\, n^{-(1-s)}\,\frac{\zeta(s)\zeta(1-s)}{s^q(1-s)^q}\, ds\, .\label{Gmdef}
\eeqa

Concerning the Dirichlet polynomials $\sumnl_{n\le N} a_n\, n^{-s}$, Bettin, Conrey, and Farmer \cite{BCF} have shown that the (array of) coefficients $\lambda_{n,N} = \mu_n (1-\log n/\log N),\ 1\le n \le N$, $\mu$ the M\"obius function, enjoy an optimality property. If RH holds and the sums $\sum_{|\Im \rho| \le T} |\zeta'(\rho)|^{-2}$ involving the critical zeros $\rho$ of $\zeta(s)$ satisfy a certain condition, then these coefficients are optimal in the sense that the related approximation error achieves the (unconditional) lower bound established by B\'aez--Duarte et al.~\cite{BBLS0} and later refined by Burnol \cite{Bu2}. With these coefficients the limit of the mixed term in (\ref{qdec}) is readily evaluated.

\begin{prop}\label{mxlim} 
$\quad \limnl_{N \to \infty} \ \sumnl_{n\le N}  \lambda_{n,N}\, F^{(q)}_n = q \quad (q=1,2).$
\end{prop}

Consequently, for both $q=1$ and $q=2$, RH is equivalent to the condition that 
the quadratic forms $\sumnl_{m,n\le N}\, \lambda_{m,N}\,  \lambda_{n,N}\,G_{m,n}^{(q)}$ 
converge to $q$ as $N\to \infty$. This assigns special interest to the Gram matrices $G_{m,n}^{(q)}$, and we will focus on these objects in the following. Finally, in Section \ref{qfasy} we derive decompositions of the associated quadratic forms  and conclude with a tentative discussion of their components' asymptotic behavior, to be illustrated in an appendix by numerical computations.

\section{Basic evaluation of Gram kernels}\label{basic}

We will consider, slightly more generally, {\em Gram kernels} $G_{u,v}^{(q)}$ defined by
\beqa
G_{u,v}^{(q)} \ceq \frac{1}{2\pi i}\int_{\Re(s) = 1/2}\, u^{-s}\, v^{-(1-s)}\,\frac{\zeta(s)\zeta(1-s)}{s^q(1-s)^q}\, ds\qquad (u,v > 0).\label{Gdef}
\eeqa
To state our evaluation of the integral (\ref{Gdef}) we introduce some notation. Let $H(x) = \sumnl_{1\le k \le x} 1/k \ \,  (x > 0)$ denote the harmonic sum function, and $\{x\}= x - \lfloor x \rfloor$ the fractional part of $x$, where $\lfloor x \rfloor$ is the largest integer $\le x$. Moreover, $\gamma_0 \equiv \gamma,\, \gamma_1$ denote the first two Stieltjes constants appearing in the Laurent series expansion of $1/\zeta(s)$. The convergence and continuity in $x$ of the functions $S_1(x),\, S_2(x)$ appearing below, here called `M\"untz type series' (e.g., \cite[Ch 2]{T}), is guaranteed by Proposition \ref{RS}.

\smallskip\noindent
\begin{thm}\label{mainthm} 
CASE $q = 1$. Let $K = (\log 2\pi - \gamma + 1)/2$, and put $S_1(x) = \sumnl_{k \ge 1} R_1(kx)$ where
\beq\label{R1def}
R_1(x) = \log x + \gamma - H(x) - \frac{\{x\}- 1/2}{x} \quad (x > 0).
\eeq
Then
\beq\label{G1eval}
G^{(1)}_{u,v} = \frac{1}{v} \left(K + \frac12 \log\frac{v}{u} \right) + \frac{1}{u}\, S_1(v/u) \qquad ( u,v >0).
\eeq
CASE $q = 2$. Let 
$$
K_1 = (\log 2\pi - \gamma + 2)/2,\quad  K_2 = \left(1-\frac{\gamma}{2}\right) \log 2\pi + \frac{1}{4} \log^2 2\pi + \frac{\pi^2}{48} - \frac{\gamma^2}{4} -\gamma - \gamma_1  + \frac{3}{2}
$$
and put $S_2(x) = \sumnl_{k \ge 1} R_2(kx)$ where
\beqa\label{R2def}
 R_2(x) \ceq \frac{1}{x} \left(\frac12 \log 2\pi + 1 + \frac12 \log x\right) + (2-\gamma)\, \log x - \frac12 \log^2 x + 2\gamma + \gamma_1 - 3 \\ 
&& \!\!\! \!\!\!\!  + \frac{1}{x}\, \sumnl_{\ell\le x} \left[\left(1 + \frac{x}{\ell}\right)\log \frac{x}{\ell} +2 \left(1 - \frac{x}{\ell}\right) \right] \quad (x > 0).\nonumber
\eeqa
Then
\beq\label{G2eval}
G^{(2)}_{u,v} = \frac{1}{v} \left(K_2 + K_1 \, \log\frac{v}{u} + \frac{1}{4} \log^2 \frac{v}{u} \right) + \frac{1}{u}\, S_2(v/u) \qquad ( u,v > 0).
\eeq
\end{thm}

\smallskip
The expression (\ref{G1eval}) is equivalent to a formula given by B\'aez--Duarte, Balazard, Landreau, and Saias \cite{BBLS} in their Proposition 90. However, our proof using straightforward residue calculus is more simple. The connection between the two expressions is detailed in Section \ref{other}. 

\smallskip
For the proof we need some facts about the functions $R_q$ and $S_q$  to be established in Section \ref{ancprf}, and about relevant residues.

\begin{prop}\label{RS}
The functions $R_q\, ,\, q=1,2$ are continuous in the range $x>0$. As $x \to \infty$,
\beqa\label{r1xp}
R_1(x) \ceq \frac{B_2(\{x\})}{2x^2} + \frac{B_3(\{x\})}{3x^3}  + O(x^{-4})\, ,\\
R_2(x) \ceq -\frac{B_4(\{x\})}{24x^4}  + O(x^{-5	})
\eeqa
where $B_2(t) = t^2 - t + 1/6,\ B_3(t) = t^3 - 3t^2/2 + t/2,\ B_4(t) = t^4 - 2t^3 + t^2 - 1/30$ are Bernoulli polynomials. \\
The series $S_q(x)$ converge absolutely for every $x>0$, are continuous in $x$, and as $x \to \infty$ we have $S_1(x) = O(x^{-2})$ and $S_2(x) =  O(x^{-4})$. 
\end{prop}

\begin{prop}\label{resids}
Let $r>0$, and for $q = 1,2$ put
\beqas 
f(s) \ceq r^s\, \frac{\zeta(s) \zeta(1-s)}{s^q(1-s)^q}\, ,\qquad
g(s) = r^s\, \frac{\zeta(s)}{s^q(1-s)^q}\, , \\ 
h(s) \ceq r^s\, \frac{\zeta(1-s)}{s^q(1-s)^q}\, ,\qquad
k(s) = r^s\, \frac{1}{s^q(1-s)^q} \, .
\eeqas 
Let $K,K_1,K_2$ be as defined above.
The residues of the functions $f,g,h,k$ at $s=0,1$ are the following.\\
CASE $q=1$. 
\beqas
\R_f(0) \ceq K + \frac12\, \log r   \qquad\qquad\quad\, \R_f(1) = -r \left(K - \frac12\, \log r\right)\\ 
\R_g(0) \ceq -\frac12  \qquad\qquad\qquad\qquad\quad \R_g(1) = r \left(1 - \gamma - \log r\right)\\
\R_h(0) \ceq -\left(1 - \gamma + \log r \right)  \qquad\quad\, \R_h(1) = \frac{r}{2}\\
\R_k(0) \ceq 1  \qquad\qquad\qquad\qquad\quad\quad\, \R_k(1) = -r
\eeqas
CASE $q=2$. 
\beqas
\R_f(0) \ceq K_2 + K_1\, \log r + \frac{1}{4}\, \log^2 r   \qquad\quad\, \R_f(1) = -r \left(K_2 - K_1\, \log r + \frac{1}{4}\, \log^2 r\right)\\
\R_g(0) \ceq -\frac12\,\left(\log 2\pi + 2 + \log r\right)  \qquad\qquad \R_g(1) = r \left(3-2\gamma - \gamma_1 + (\gamma - 2)\log r  + \frac12\, \log^2 r\right)\\
\R_h(0) \ceq -3+2\gamma + \gamma_1 + (\gamma - 2)\log r  - \frac12\, \log^2 r \qquad\quad\! \R_h(1) = \frac{r}{2} \left(\log 2\pi + 2 - \log r\right)\\
\R_k(0) \ceq 2+\log r  \qquad\qquad\qquad\qquad\quad\qquad\! \R_k(1) = -r\left(2-\log r\right)
\eeqas
\end{prop}

\medskip\noindent
{\bf {\em Proof of Theorem \ref{mainthm}. }}  
Throughout the following, integration along vertical lines $\Re(s) = c$ is understood to go upwards from $c -i\infty$ to $c+i\infty$ and is simply denoted as $\int_c$. It is known \cite[Sect. 9.2]{Ed} that for fixed $\sigma \in \real$ and every $\el >0$ one has for $|t|\to \infty$
\beqas
\zeta(\sigma+it) \ceq \left\{
\begin{array}{ll} 
O(|t|^{1/2-\sigma}) & \ \mbox{if}\quad \sigma < 0 \\ 
O(|t|^{(1-\sigma)/2 + \el}) & \ \mbox{if}\quad 0 \le \sigma \le 1 \\  O(1) & \ \mbox{if}\quad \sigma > 1.
\end{array} \right.
\eeqas
Consequently, the integral (\ref{Gdef}) stays absolutely convergent when the line of integration is shifted to any line $\Re(s) = c\notin \{0,1\}$ such that, respectively, $-1/2 < c < 3/2$ if $q=1$, or $-5/2 < c < 7/2$ if $q=2$. Furthermore, the horizontal contributions to the contour integral at large imaginary values are negligible.  

We now consider the case $q=1$. Shifting the line of integration at first to $\Re(s) = -1/4 \equiv c_-$, say, we obtain 
\beq\label{gram1}
G^{(1)}_{u,v} = \frac{1}{v}\left(K + \frac12 \log(v/u)\right) + J(c_-),
\eeq
where the first term is $\mbox{Res}\left(u^{-s} v^{-(1-s)}\,\frac{\zeta(s)\zeta(1-s)}{s(1-s)}\right) \bigg|_{s=0}$, the residue of the integrand at $s=0$, and
\beqa
J(c_-) \ceq \frac{1}{2\pi i}\,\int_{c_-}\, u^{-s} v^{-(1-s)}\, \frac{\zeta(s)\zeta(1-s)}{s(1-s)}\, ds\nonumber\\ \ceq
\frac{1}{2\pi i}\,\int_{c_-}\, \sumnl_{k\ge 1}\, u^{-s} (kv)^{-(1-s)}\, \frac{\zeta(s)}{s(1-s)}\, ds\nonumber\\ \ceq
\sumnl_{k\ge 1}\, \frac{1}{kv}\,J_k(c_-), \label{critxch}
\eeqa
with
$$
J_k(c_-) = \frac{1}{2\pi i}\,\int_{c_-}\, (kv/u)^s\, \frac{\zeta(s)}{s(1-s)}\, ds.
$$
The interchange of summation and integration at (\ref{critxch}) is allowed because the inner sum is absolutely convergent and integrable, and the last sum converges, too; see below.\\
We next shift to the line $\Re(s) = 5/4 \equiv c_+$. Picking up (minus) the residues of $(kv/u)^s\, \frac{\zeta(s)}{s(1-s)}$ at zero and one, then expanding $\zeta(s)$ we similarly get
\beq\label{xch2}
J_k(c_-) = \frac12 - \frac{kv}{u}\left(1- \log\frac{kv}{u} - \gamma\right) + \sumnl_{\ell \ge 1}\, J_{k,l}(c_+)
\eeq
where 
\beqas
J_{k,l}(c_+) \ceq \frac{1}{2\pi i}\,\int_{c_+}\, \left(\frac{kv}{\ell u}\right)^{\! s}\, \frac{ds}{s(1-s)}\, .
\eeqas
Here sum and integral may be interchanged because the sum in (\ref{xch2}) is finite: $J_{k,l}(c_+)$ vanishes if $kv/(\ell u) < 1$, or $\ell > kv/u$. For $\ell \le kv/u$ it equals $1- \frac{kv}{\ell u}$. \\
Collecting terms one obtains 
\beqas
\frac{1}{kv}\,J_k(c_-) \ceq 
\frac{1}{kv}\, \left[\, \frac12 - \frac{kv}{u}\left(1- \log\frac{kv}{u} - \gamma\right) +\!\!\! \sum_{1 \le \ell \le kv/u} \left(1- \frac{kv}{\ell u}\right)\right]\\
\ceq \frac{1}{u}\, \left[\, \frac{u}{2kv} - \left(1- \log\frac{kv}{u} - \gamma\right) + \bigg\lfloor\frac{kv}{u}\bigg\rfloor\, \frac{u}{kv} -  H(kv/u)\, \right] \\ \ceq
\frac{1}{u}\, \left[\, \log\frac{kv}{u} + \gamma - H(kv/u) + \left(\frac12 - \left\{\frac{kv}{u}\right\}\right)\bigg/ \frac{kv}{u}\, \right] \\ \ceq
\frac{1}{u}\, R_1(kv/u) .
\eeqas
Now by Proposition \ref{RS} the sum $\sumnl_{k\ge 1} \frac{1}{kv}\,J_k(c_-)$ is convergent, which settles the interchange at (\ref{critxch}). But $J(c_-) = u^{-1} \sumnl_{k \ge 1} R_1(kv/u) = u^{-1}\,S_1(v/u)$, so (\ref{G1eval}) follows from (\ref{gram1}).

The case $q=2$ can be treated along exactly the same lines. Only the respective residues differ. \done

\section{Further representations, reciprocity relations}\label{furtherrep}

The following proposition points out that the kernel $G^{(2)}$ is in fact a scale average of the kernels $G_{u,v}^{(1)}$. The integral representation of the latter is well-known (e.g., \cite[p.~38]{BBLS}, \cite[p.~5714]{BC}).

\begin{prop}\label{intrep}
For any $u,\, v>0$
\beqa\label{CMG1}
G_{u,v}^{(1)} \ceq \frac{1}{uv} \int_0^\infty \{tu\} \{tv\}\, \frac{dt}{t^2}\, .\\
G_{u,v}^{(2)} \ceq \int_0^1 \int_0^1 G_{ux,vy}^{(1)}\, dx\, dy \, .\label{CMG2} 
\eeqa
\end{prop}
{\bf {\em Proof.} } 
We follow Conrey and Myerson \cite{CM}, making use of some facts about Mellin transforms. Putting technicalities aside, if $\widehat f(s) = \int_0^\infty f(x)\, x^{s-1}\, dx$ denotes the Mellin transform of a function $f$ on $(0,\infty)$, then the Mellin transform of the function $\I f(x) = f(1/x)/x$ is $\widehat f(1-s)$, and the Mellin transform of the (multiplicative) convolution $f \ast \I f$ is $\widehat f(s)\, \widehat f(1-s)$. In particular, if $f(x) = \{1/x\}$, which has the Mellin transform
\beq\label{mellfrc}
\widehat f(s) = \int_0^\infty \left\{\frac{1}{x}\right\} \, x^{s-1}\, dx = -\frac{\zeta(s)}{s}\, ,
\eeq
then Mellin inversion gives
$$
G_{u,v}^{(1)} = \frac{1}{v}\, (f \ast \I f)(u/v) = \frac{1}{v} \int_0^\infty  \frac{\{x\}}{x}\,\left\{\frac{x}{u/v}\right\} \,\frac{dx}{x} = \frac{1}{uv} \int_0^\infty \{tu\} \{tv\}\, \frac{dt}{t^2}\, .
$$
As for $q=2$, note first that the Mellin transform of the indicator function $\chi$ of the unit interval is $\widehat \chi(s) = 1/s$, so that the Mellin transform of the convolution $h = \chi\ast\{1/\cdot\}$ is $-\zeta(s)/s^2$.
Thus, Mellin inversion as above gives, after some calculation,
$$
G_{u,v}^{(2)} = \frac{1}{v}\, (h\ast \I h) (u/v) = \frac{1}{uv} \int_0^\infty \int_0^1 \{tux\}\, \frac{dx}{x} \int_0^1 \{tvy\}\, \frac{dy}{y}\, \frac{dt}{t^2}\, . 
$$
By Fubini this may also be written as 
$$
G_{u,v}^{(2)} = \int_0^1 \int_0^1 \int_0^\infty \frac{\{tux\}}{ux}\, \frac{\{tvy\}}{vy}\, \frac{dt}{t^2}\, dx\, dy = \int_0^1 \int_0^1 G_{ux,vy}^{(1)}\, dx\, dy \, ,
$$
which completes the proof. \done

\newpage
Incidentally, (\ref{CMG1}), (\ref{CMG2}) show that the kernels $G_{u,v}^{(q)}$ are symmetric in $u,v$, which is not obvious from Theorem \ref{mainthm}. Manifestly symmetric expressions can be obtained by convex combination of $G^{(q)}$ and its transpose. 

\begin{prop}\label{convcomb}
\beqa\label{G1cc}
G^{(1)}_{u,v}\ceq \frac{1}{u+v}\, \left[\,2K + \frac{v}{u}\, S_1(v/u) + \frac{u}{v}\, S_1(u/v)\, \right],\\
G^{(2)}_{u,v} \ceq \frac{1}{u+v}\, \left[\,2K_2 + \frac12\, \log^2\frac{v}{u} + \frac{v}{u}\, S_2(v/u) + \frac{u}{v}\, S_2(u/v)\, \right] ,\label{G2cc}\\
G^{(1)}_{u,v} \ceq \frac12 \left[\, K \left(\frac{1}{v} + \frac{1}{u}\right) + \frac12 \left(\frac{1}{v} - \frac{1}{u}\right) \log \frac{v}{u}  + \frac{S_1(u/v)}{v} +  \frac{S_1(v/u)}{u}\, \right], \label{G1eqw}\\ 
G^{(2)}_{u,v} \ceq \frac12 \left[\left(K_2+\frac{1}{4}\, \log^2 (v/u)\right) \left(\frac{1}{v} + \frac{1}{u}\right) + K_1 \left(\frac{1}{v} - \frac{1}{u}\right) \log \frac{v}{u}  + \frac{S_2(u/v)}{v} +  \frac{S_2(v/u)}{u}\, \right]  . \label{G2eqw}
\eeqa
\end{prop}
{\bf {\em Proof.} } 
The particular choice $\frac{v}{u+v}G^{(q)}_{u,v}+ \frac{u}{u+v}G^{(q)}_{v,u}$ deletes the respective second terms in (\ref{G1eval}), (\ref{G2eval}) and gives the first two expressions. The second pair is obtained by taking equal weights $1/2$ each. \done

\smallskip
By equating (\ref{G1eval}) and (\ref{G1cc}) (or (\ref{G1eqw})) and setting $r=v/u$, it is possible to represent $S_1(1/r)$ in terms of $S_1(r)$ and known quantities, and analogously for $S_2(1/r)$. This readily yields the following reciprocity formulae.

\begin{cor}\label{reci}
{\em [Reciprocity relations]} For every $r>0$
\beqa\label{reci1}
S_1(1/r) \ceq rS_1(r) + K\, (1-r) + \frac12\, (1+r)\, \log r \, ,\\
S_2(1/r) \ceq rS_2(r)  + K_2\, (1-r) + K_1\, (1+r)\, \log r + \frac{1}{4}\, (1-r)\,  \log^2 r\, . \label{reci2}
\eeqa
\end{cor}

Putting these findings together one can obtain yet another representation of the kernel $G^{(2)}$ akin to (\ref{G1eqw}) which we state here without proof.

\smallskip\noindent
\begin{thm}\label{G2thm}
With $K_1 = K+1/2$ as defined in Theorem \ref{mainthm} we have
\beqa
G^{(2)}_{u,v} \ceq 
K_1 \left(\frac{1}{v} + \frac{1}{u}\right) + \frac12 \left(\frac{1}{v} - \frac{1}{u}\right) \log \frac{v}{u} + \frac{1}{v} \int_{u/v}^\infty \frac{S_1(r)}{r}\, dr + \frac{1}{u} \int_{v/u}^\infty \frac{S_1(r)}{r}\, dr  \, .\label{G2alt2}
\eeqa
\end{thm}

\section{Comparison with representations from the literature}\label{other}

The case $q=2$ does not seem to have been considered before, so we focus on $q=1$. 
As noted earlier the expression (\ref{G1eval}) for $G^{(1)}_{m,n}$ appears in a different form already in \cite[Proposition 90]{BBLS}. It reads 
\beq\label{bbls90}
A(r) = K + \frac12 \log r - r \int_r^\infty \varphi_1(t)\frac{dt}{t^2} \quad (r > 0).
\eeq
Here $A(r) = \int_0^\infty \{xr\} \{x\}\, \frac{dx}{x^2},\, r >0$ denotes the `autocorrelation' of the fractional parts functions $\{x\} = x - \lfloor x \rfloor$ and $\{xr\}$, and  $\varphi_1(t),\, t >0$ is defined almost everywhere as 
\beq\label{fi1def}
\varphi_1(t)  = \sumnl_{k\ge 1} (\{kt\} - 1/2)/k\, .
\eeq
Since by (\ref{CMG1}) one has $A(n/m)/n = G^{(1)}_{m,n}$, we may write (\ref{bbls90}) as 
$$
G^{(1)}_{m,n} = \frac{1}{n} \left(K + \frac12 \log\frac{n}{m} \right) - \frac{1}{m}\, \int_{n/m}^\infty \varphi_1(t)\, \frac{dt}{t^2}\, ,
$$
whence by comparison with (\ref{G1eval}) it must hold that 
\beq\label{S1rel}
S_1(r) =-\int_r^\infty \varphi_1(t)\, \frac{dt}{t^2}\, , \quad r = n/m.
\eeq
For a direct proof of (\ref{S1rel}) valid for all $r>0$ see Proposition \ref{S1rep1}. 

\smallskip
Let us point out that {\em given (\ref{S1rel}), the reciprocity relation for  $S_1$ is also immediate from (\ref{bbls90}) via the corresponding reciprocity relation $A(r)=r A(1/r)$ for the autocorrelation function} $A$ (which latter is just a change of variables).

\smallskip
In \cite[Proposition 88]{BBLS} and \cite[Proposition 89]{BBLS} the authors also gave alternative representations of the autocorrelation function $A$, namely 
$$
A(r) = \frac{1-r}{2}\, \log r + \frac{1+r}{2}\, (\log 2\pi - \gamma)  -  \varphi_1(r) - r \varphi_1(1/r),
$$
valid for such $r>0$ where the series $\varphi_1(r)$ converges, and 
\beq\label{bbls89}
A(r) = \frac{1-r}{2}\, \log r + \frac{1+r}{2}\, (\log 2\pi - \gamma)  -  \frac{\pi}{2m} \left(V(n,m)+ V(m,n)\right),
\eeq
valid for rational $r = n/m$ such that $n$ and $m$ have no common divisor. Here 
$$
V(n,m) = \sum_{k=1}^{m-1} \left\{\frac{kn}{m}\right\} \cot\left(\frac{\pi k}{m}\right)
$$
is a `Vasyunin sum'. Such sums appeared at first in Vasyunin's \cite{V} evaluation of a closely related Gram matrix, and are now being studied on their own right due to their reciprocity properties and connections with Eisenstein series; see e.g. \cite{BC,LZ}.  
In passing, the relations (\ref{bbls90}), (\ref{bbls89}) along with (\ref{S1rel}) can be used to show that the Vasyunin term $-\pi (V(n,m) + V(m,n))/(mn)$ can be written in terms of our M\"untz type series as $S_1(n/m)/m + S_1(m/n)/n$ plus known functions. 

B\'aez--Duarte et al.~also determined the exact value of $A(1) = G^{(1)}_{1,1}$, which is $G^{(1)}_{1,1} = \log 2\pi - \gamma$; see \cite[Proposition 87]{BBLS}. Since $G^{(1)}_{1,1} = K + S_1(1)$ and $K = (\log 2\pi - \gamma + 1)/2$ it follows that $S_1(1) = (\log 2\pi - \gamma - 1)/2$. For an independent evaluation of $S_1(1)$ see Proposition \ref{S11}.

\section{Integral representations of the series $S_q$}\label{Sreps}

We begin by giving an independent proof of the representation of $S_1$ obtained indirectly in Section \ref{other}.

\begin{prop}\label{S1rep1}\ 
\beq\label{S1rep}
S_1(r) = -\int_r^\infty \varphi_1(t)\, \frac{dt}{t^2}\, , \quad r >0.
\eeq
\end{prop}
{\bf {\em Proof.} } 
We first show that
\beq\label{R1rep}
R_1(x) = -\int_x^\infty (\{t\}-1/2)\, \frac{dt}{t^2}\, , \quad x >0.
\eeq
By its definition, $R_1 = f+g$ where
$$
f(x) = \frac{1}{2x} + \log x + \gamma - 1,\qquad g(x) = -H(x) + \frac{\flx}{x}\, .
$$
Clearly $g(x)=0$ for $x \le 1$, while for $1\le n < x < n+1$ we have $g(x) = -H(n) + n/x$, hence $g'(x) = -n/x^2 = -\flx / x^2$. Now $f'(x) = -\frac{1}{2x^2} + \frac{1}{x}$, so 
$
R_1'(x) = -\frac{1}{2x^2} + \frac{1}{x}-\frac{\flx}{x^2} = \frac{x- \flx - 1/2}{x^2}\, ,
$
and the representation (\ref{R1rep}) follows by integration. Indeed, $R_1$ is continuous (Proposition \ref{RS}); and a possible integration constant must vanish because both sides of (\ref{R1rep}) tend to zero as $x \to\infty$. 
The representation (\ref{S1rep}) then follows from 
\beq\label{SSrep}
S_1(r) = -\sumnl_{k\ge 1} \int_{kr}^\infty (\{t\}-1/2)\, \frac{dt}{t^2} = 
-\int_{r}^\infty \sumnl_{k\ge 1} \frac{\{kx\}-1/2}{k}\, \frac{dx}{x^2}\, . 
\eeq
It remains to justify the interchange of summation and integration. For differentiable $h$ and any $x$ we have 
$$
\left| \int_x^{x+1} (\{t\}-1/2)\, h(t)\, dt\,  \right| \le  \frac12 \supnl_{x \le t \le x+1} |h'(t)|\, .
$$
Thus if $m$ denotes the smallest integer $>= kr$ we have 
\beqas
\left| \int_{kr}^\infty (\{t\}-1/2)\, \frac{dt}{t^2} \right| \cleq  \frac{1}{2(kr)^2} + \sumnl_{n\ge m} \left| \int_{n}^{n+1} (\{t\}-1/2)\, \frac{dt}{t^2}\right| \le \frac{1}{2(kr)^2} +   \sumnl_{n\ge m} n^{-3} \\ \ceq O( (kr)^{-2})\, ,
\eeqas 
which establishes the absolute convergence of the first series in (\ref{SSrep}). 
Each term in the last series, denoted $\varphi_1(x)$, is periodic with period 1, hence so is $\varphi_1(x)$. The formal Fourier series expansion is
\beqas
\varphi_1(x) \ceq \sumnl_{k\ge 1} \frac{\{kx\}-1/2}{k} = -\sumnl_{k \ge 1}\sumnl_{n \ge 1} \frac{\sin 2\pi nkx}{ \pi nk} = - \sumnl_{m \ge 1} \frac{d(m)}{\pi m}\, \sin 2\pi mx
\eeqas 
where $d(m)$ denotes the divisor function (which counts all pairs of natural numbers whose product equals $m$). Now $\sum_{m \ge 1} \frac{d(m)^2}{ m^2}$ is convergent because $d(m) = O(m^\el)$ for any fixed $\el > 0$ \cite[Theorem 315]{HW}. Therefore $\varphi_1 \in L^2(0,1)$, and the existence of the last integral in (\ref{SSrep}) follows by Cauchy-Schwarz,
\beqas
\int_1^\infty |\varphi_1(x)|\, \frac{dx}{x^2} \ceq \sumnl_{n \ge 1} \int_n^{n+1} |\varphi_1(x)|\, \frac{dx}{x^2} \le 
\sumnl_{n \ge 1}  \left( \int_0^1 \varphi_1(x)^2\, dx \, \int_n^{n+1}x^{-4}\, dx\right)^{1/2} < \infty\, ;
\eeqas
clearly also, $\int_r^1 |\varphi_1(x)|\, \frac{dx}{x^2} < \infty$ for every $0 <r < 1$. The last part draws on \cite[Section 7]{BBLS}. \done

\medskip
The case $q=2$ can be treated similarly.

\begin{prop}\label{S2rep}
Let 
\beq\label{Vdef}
V(x) = x\, \{H(x) - \log x - \gamma + 2\} - \frac12\log x\,  +\, \sumnl_{n \le x}\! \log n\, - \flx \, (1+\log x) -\frac12 \log 2\pi - \frac12\, .
\eeq
Then for $r >0$,
\vspace*{-2mm}
\beqa\label{S2int}
S_2(r) \ceq \int_{r}^\infty S_1(t) \left( \frac{1}{t} -  \frac{1}{r}\right) dt\, ,\\ 
\int_{r}^\infty S_1(t)\, dt\ceq \sumnl_{n \ge 1} \frac{V(nr)}{n} \, . \label{VS1}
\eeqa
\end{prop}
{\bf {\em Proof. }}
The formula 
\beq\label{inteval}
\int_1^\infty \left(\{x\} - 1/2\right) \frac{dx}{x} = \frac12 \log 2\pi - 1,
\eeq
needed below is certainly known. A proof can be based on the well-known relation \cite[p.~14]{T}
\beqas
-\frac{\zeta(s)}{s} \ceq 
\frac{1}{1-s} + \frac{1}{2s} + \int_1^\infty \left(\{x\}-1/2\right) x^{-1-s}\, dx\qquad (\Re(s) > -1).\label{classical}
\eeqas
Expanding $\frac{\zeta(s)}{s}+\frac{1}{1-s} + \frac{1}{2s}$ at $s=0$, then letting $s$ tend to zero gives (\ref{inteval}). We now will show that
\beqa\label{Vrep}
V(x) \ceq \int_x^\infty R_1(t) \, dt\, . 
\eeqa
The function $V(x)$ from (\ref{Vdef}) is continuous, and for $x$ not an integer its derivative equals
\beqas
V'(x) \ceq H(x) - \log x - \gamma + 2 - 1 -\frac{1}{2x} - \frac{\flx}{x} = H(x) - \log x - \gamma + \frac{x - \flx - 1/2}{x} \nonumber\\ \ceq -R_1(x) \, \label{V1rep}
\eeqas
which is continuous itself. Therefore $V(x) = V(1) -\int_1^\infty R_1(t)\, dt + \int_x^\infty R_1(t)\, dt$. Now by Fubini
\beqas
\int_1^\infty R_1(t)\, dt  \ceq - \int_1^\infty \int_t^\infty (\{u\}- 1/2)\, \frac{du}{u^2}\, dt = - \int_1^\infty \left(\int_1^u \, dt\right) (\{u\}- 1/2)\, \frac{du}{u^2}\\ \ceq
- \int_1^\infty (u-1) \, (\{u\}- 1/2)\, \frac{du}{u^2} = \int_1^\infty (\{u\}- 1/2)\, \frac{du}{u^2} -\int_1^\infty (\{u\}- 1/2)\, \frac{du}{u} \\ \ceq
- R_1(1) - \left(\frac12\log 2\pi - 1\right) = \frac12 - \gamma - \frac12\log 2\pi + 1 = \frac{3}{2} - \gamma - \frac12\log 2\pi \\ \ceq V(1)\, ,
\eeqas
wherein we have used (\ref{inteval}). Thus, the integration constant $V(1) -\int_1^\infty R_1(t)\, dt$ vanishes, and (\ref{Vrep}) is established. The relation (\ref{VS1}) then follows by an (admissible) interchange of summation and integration, taking into account the definition of $S_1$.

\noindent
Differentiation of $R_2(x)$ as defined in (\ref{R2def}) gives $R_2\,\!\!\! '(x) = V(x)/x^2$, whence on arguing as previously,
\beqa\label{R2rep2}
R_2(x) \ceq -\int_x^\infty V(t)\, \frac{dt}{t^2} = \int_x^\infty V(t)\, \left(\frac{1}{t}\right)'\, dt 
= \int_x^\infty R_1(t)\, \frac{dt}{t} -\frac{V(x)}{x}\, .
\eeqa
Using this relation along with the definitions of $S_2$ and $S_1$ we obtain
\beqa
S_2(r) + \sumnl_{n \ge 1} \frac{V(nr)}{nr} \ceq  \sumnl_{n \ge 1} \int_{nr}^\infty R_1(x)\, \frac{dx}{x}  = \sumnl_{n \ge 1} \int_{r}^\infty R_1(nt)\, \frac{dt}{t} =  \int_{r}^\infty S_1(t)\, \frac{dt}{t}\, \label{S2frst}
\eeqa
which in view of (\ref{VS1}) implies (\ref{S2int}). \done


\medskip
We conclude this section by stating yet another representation of $S_1(r)$ without proof.

\begin{prop}\label{S1rep2}
\beq\label{SCM}
S_1(r) = -\frac{1}{r} \int_1^\infty \left(\{xr\} -\frac12\right) \lfloor x \rfloor\, \frac{dx}{x^2}\, , \quad r> 0.
\eeq
\end{prop}

\medskip
Together with Proposition \ref{S1rep1} this yields the peculiar identity

\beqa\label{syn1}
\int_1^\infty \left(\sumnl_{k\ge 1} \frac{\{kxr\} - \frac12}{k}\right) \frac{dx}{x^2} \ceq  \int_1^\infty \left(\{xr\} -\frac12\right) \lfloor x \rfloor\, \frac{dx}{x^2}\, . 
\eeqa

\section{The constants $S_q(1)$}\label{Sq1}

In view of (\ref{G1eval}), (\ref{G2eval}), an exact evaluation of $G^{(q)}_{1,1}$ is possible if $S_q(1)$ can be calculated explicitly. For $q=1$ this was achieved in \cite[Proposition 87]{BBLS}. Here is an independent evaluation.

\begin{prop}\label{S11}
\beq\label{S1eval}
S_1(1) = (\log 2\pi - \gamma - 1)/2.
\eeq
\end{prop}
{\bf {\em Proof. }}
Recall that $S_1(x) = \sumnl_{n \ge 1} R_1(nx)$ where
$$
R_1(x) = \log x + \gamma - H(x) - \frac{x - \lfloor x \rfloor - 1/2}{x} \quad (x > 0).
$$
We have
\beqas
\sumnl_{n \le N}\, H(n) \ceq \sumnl_{n \le N}\, \sumnl_{k \le n}\, \frac{1}{k}\, =\, \sumnl_{k \le N}\,  \frac{1}{k}\, \sumnl_{k \le n \le N} 1\, =\,  \sumnl_{k \le N}\,  \frac{1}{k}\, (N-k + 1) \nonumber\\ \ceq (N+1)\, H(N) - N\, .\label{H0sum}
\eeqas
Using the well-known expansion for the harmonic series and Stirlings's formula,
\beqas
H(N) \ceq \log N + \gamma + \frac{1}{2N} + O(N^{-2})\, ,\\
\sumnl_{n\le N} \log n  \ceq N\log N - N + \frac12 \log N + \frac12 \log 2\pi + O(N^{-1})\, ,
\eeqas
we get
\beqas
\sumnl_{n\le N} R_1(n) \ceq \sumnl_{n\le N} \log n  + N\gamma - (N+1)\, H(N) + N + \frac12\, H(N) \\ \ceq
N\log N - N + \frac12 \log N + \frac12 \log 2\pi + O(N^{-1}) +  N\gamma +N \\ \!\!\! && -\, (N+\frac12) \,  \left(\log N + \gamma + \frac{1}{2N} + O(N^{-2})\right)\\ \ceq
\frac12 \log 2\pi - \frac12 \gamma - \frac12 + O(N^{-1})\, .
\eeqas
Letting $N$ tend to infinity completes the proof. \done

\medskip
Unfortunately, this simple approach does not work when $q=2$. At least we can state the following.

\begin{prop}\label{S21}
\beqa\label{S21rep}
S_2(1) \ceq \int_1^\infty \frac{S_1(r)}{r}\, dr - \int_1^\infty S_1(r)\, dr = K_2 - 2K_1 - 2\int_1^\infty S_1(r)\, dr\, .
\eeqa
\end{prop}
{\bf {\em Proof. }}
The first relation is immediate from (\ref{S2int}). Putting $u=v=1$ in (\ref{G2cc}) and (\ref{G2alt2}), respectively, we obtain
$$
G_{1,1}^{(2)} = K_2 + S_2(1) = 2 K_1 + 2 \int_1^\infty \frac{S_1(r)}{r}\, dr \, ,
$$
whence the second relation follows by equating the two expressions for $S_2(1)$. \done 

\smallskip
A numerical evaluation based on the rapidly converging series $\sum_{n\ge 1} V(n)/n$ (which by (\ref{VS1}) is equal to $\int_1^\infty S_1(r)\, dr$)  gives the approximative values
$$
S_2(1) \doteq 0.000643, \quad G_{1,1}^{(2)} \doteq 3.270465.
$$
In the case $q=1$,
$$
S_1(1) \doteq 0.130331, \quad G_{1,1}^{(1)} \doteq 1.260661.
$$

\section{Ancillary proofs}\label{ancprf}

\subsection{Proof of Proposition \ref{mxlim}}

Integrating over large half circles to the right of the line $\Re s =1$ one finds that $F^{(q)}_n$ equals minus the residue of $n^{-s}\, \zeta(s) / (s(1-s))^q$ at $s=1$, so that by Proposition \ref{resids}
$$
F^{(1)}_n = \frac{1}{n} \left(\gamma -1 - \log n\right), \quad F^{(2)}_n = \frac{1}{n} \left(-\frac12\log^2 n + (\gamma - 2)\log n + 2 \gamma +\gamma_1 - 3\right).
$$
Furthermore, the limits of the summatory functions $L_{\mu,k}(x) = \sumnl_{n\le x}  \frac{\mu(n)}{n}\, \log^k n,\ k \le 2$ are known: one has 
\beqa
\limnl_{x\to\infty} L_{\mu,0}(x) \ceq 0\, ,\label{L0lim}\\
\limnl_{x \to\infty} L_{\mu,1}(x) \ceq -1\, ,\label{L1lim}\\
\limnl_{x \to\infty} L_{\mu,2}(x) \ceq -2\gamma\, .\label{L2lim}
\eeqa
Thus with $\lambda_{n,N} = \mu_n (1-\log n/\log N)\equiv \lambda_n$ and $L_{\mu,k} \equiv L_k$ it follows that 
\beqas
\sumnl_{n\le N} \lambda_n F^{(1)}_n \ceq \sumnl_{n\le N} \left(\gamma -1 - \log n\right) \frac{\mu_n}{n} \left(1 - \frac{\log n}{\log N}\right)\\ \ceq
(\gamma -1) L_0(N) -(\gamma -1)\, \frac{L_1(N)}{\log N} - L_1(N) +  \frac{L_2(N)}{\log N} 
\eeqas
tends to 1 as $N\to\infty$, which settles the case $q=1$. In the case $q=2$
\beqas
\sumnl_{n\le N} \lambda_n F^{(2)}_n \ceq - \frac12 \sumnl_{n\le N} \left(1 - \frac{\log n}{\log N}\right) \frac{\mu_n \log^2 n}{n} + (\gamma -2) \sumnl_{n\le N} \left(1 - \frac{\log n}{\log N}\right) \frac{\mu_n \log n}{n} \\ && \!\!\! 
+ \ (2 \gamma +\gamma_1 - 3) \sumnl_{n\le N} \left(1 - \frac{\log n}{\log N}\right) \frac{\mu_n}{n}\\ \ceq
- \frac12 \left[L_2(N)- \frac{L_3(N)}{\log N}\right] + (\gamma - 2)\left[L_1(N)- \frac{L_2(N)}{\log N}\right] + (2 \gamma +\gamma_1 - 3)\left[L_0(N)- \frac{L_1(N)}{\log N}\right]
\eeqas
converges to $-\frac12 (-2\gamma) - (\gamma -2) =2$, as claimed, provided that $L_3(N)/\log N$ tends to zero. Indeed, integration/summation by parts gives 
\beq\label{L2L3}
\frac{L_3(N)}{\log N} = L_2(N)- \frac{1}{\log N}  \int_1^N L_2(t)\, \frac{dt}{t}
 = L_2(N)+2\gamma - \frac{1}{\log N}  \int_1^N (L_2(t)+2\gamma)\, \frac{dt}{t}\, ,
\eeq
and since $L_2(N)+2\gamma$ tends to zero as $N\to\infty$, so does its logarithmic mean. %
\done

\subsection{Proof of Proposition \ref{RS}}

Let $b = \lfloor x\rfloor = x - \{x\},\ \el = \{x\}/x$. By the known expansion of the harmonic series one has for $x \to \infty$
\beqas
&& \!\!\! \log x + \gamma - H(x)\\ \ceq \log x - \log b - (2b)^{-1} + (12b^2)^{-1} + O(b^{-4}) = 
\log \frac{x}{b} - \frac{1}{2x}\,\frac{x}{b} + \frac{1}{12x^2} \left(\frac{x}{b}\right)^2 + O(x^{-4})\\ \ceq 
-\log(1-\el) - \frac{1}{2x}\,\frac{1}{1-\el} + \frac{1}{12x^2}\,\frac{1}{(1-\el)^2}+ O(x^{-4})\\ \ceq 
\el + \frac{\el^2}{2}  + \frac{\el^3}{3} -\frac{1}{2x}\left(1+ \el + \el^2\right) + \frac{1}{12x^2} \left(1 + 2\el\right) + O(x^{-4})\\ \ceq 
\frac{\{x\}-1/2}{x} + \frac{1}{2x^2}\left(\{x\}^2-\{x\}+1/6\right) + \frac{1}{6x^3}\left(2\{x\}^3-3\{x\}^2+\{x\}\right) +  O(x^{-4})\, .
\eeqas
Rearranging terms gives (\ref{r1xp}). The corresponding result for $R_2$ can be obtained by very tedious calculations using the expansions for the harmonic series and  Stirling's formula up to terms of order $O(x^{-5})$, as well as the expansion
$$
\sumnl_{k\le b}\, \frac{\log k}{k} = \frac12 \log^2 b + \gamma_1 + \frac{\log b}{2b} + \frac{1-\log b}{12b^2} + \frac{6\log b - 11}{720 b^4} + O(b^{-5})
$$
which can be derived by Euler--Maclaurin summation.\\
(Jump) discontinuities of $R_1\, ,R_2$ could occur at the natural numbers only. One readily checks  that the left-hand and the right-hand limits are identical, so the $R_q$ are continuous everywhere. Continuity of $S_1\, ,  S_2$ then ensues from the locally uniform convergence of the series. As for the latter, recall that $R_q(kx) = O((kx)^{-2q})$, which also implies the tail estimates $S_q(x) = O(x^{-2q})$ ($q = 1,2$). \done

\section{Asymptotics of the quadratic forms: Some tentative steps}\label{qfasy}

\subsection{Decomposition of the quadratic forms}\label{secqdec}

Towards a potentially useful decomposition of the quadratic forms $Q_N^{(q)} = \sumnl_{n \le N}  a_m\, a_n\, G_{m,n}^{(q)}$ based on Theorem \ref{mainthm} 
we introduce for fixed, initially arbitrary coefficients $a_n$ and $j=0,1,\ldots$
the expressions
\beqas
L_j(N) \equiv L_{a,j}(N) \ceq \sumnl_{n\le N} \frac{a_n}{n} \log^j n\, ,\\ 
M_j(N) \equiv M_{a,j}(N)\ceq \sumnl_{n\le N} a_n \log^j n\, ,\\
E^{(q)}(N) \equiv E_a^{(q)}(N) \ceq \sumnl_{m \le N} \frac{a_m}{m} \left( \sumnl_{n\le N}\, a_n\, S_q(n/m) - R_q(1/m)\right).
\eeqas
The last term is defined with a view to the M\"obius inversion formula $\sumnl_{n= 1}^\infty\, \mu_n S_q(n/m) = R_q(1/m)$. Estimation of the inversion error $E^{(q)}(x)$ is a major challenge; we set it aside. Henceforth we focus on the $L$(andau)- and $M$(ertens)-type partial sums, and write $\sumnl'_n$ to denote $\sum_{n=1}^N$. 

\medskip
\noindent
CASE $q=1$. Using (\ref{G1eval}) we may write
\beqas
Q_N^{(1)} \ceq \sumnl'_{m,n} \left(K + \frac12 \log n -\frac12 \log m \right) a_m\, \frac{a_n}{n}  + \sumnl'_{m}\, \frac{a_m}{m}\left( \sumnl'_{n}\, a_n\, S_1(n/m)\right)\nonumber\\ \ceq
\left(\sumnl'_{m} \, a_m\right) \left[\frac12\sumnl'_{n} \frac{a_n \log n}{n}  + K\sumnl'_{n}\, \frac{a_n}{n}\, \right] - \frac12 \left(\sumnl'_{m}\, a_m\, \log m\, \right) \left( \sumnl'_{n}\, \frac{a_n}{n}\right)  \nonumber\\ && \!\!\! 
+ \ \sumnl'_{m}\, \frac{a_m}{m}\, R_1(1/m) + \sumnl'_{m}\, \frac{a_m}{m} \left(\sumnl'_{n}\, a_n\, S_1(n/m) -  R_1(1/m)\right)
\nonumber\\ \ceq
M_0(N) \left(\frac12 L_1(N) + K L_0(N)\right) - \frac12 M_1(N) L_0(N) + \sumnl'_{n}\, \frac{a_m}{m}\, R_1(1/m) + E^{(1)}(N)\, .
\eeqas
Now $R_1(x) = \frac{1}{2x} + \log x + \gamma -1$ for $x \leq 1$, so
\beqas
\sumnl'_{m}\, \frac{a_m}{m}\, R_1(1/m) \ceq  \sumnl'_{m}\, \frac{a_m}{m} \left(\frac{m}{2} - \log m\right) + (\gamma - 1) \sumnl'_{m}\, \frac{a_m}{m} \\ \ceq
\frac12 M_0(N) -L_1(N) + (\gamma -1) L_0(N)\, , 
\eeqas
whence
\beq\label{Q1eval}
Q_N^{(1)} = M_0(N) \left(K L_0(N) + \frac12\, \{L_1(N) + 1\}\right) - \frac12 M_1(N)L_0(N) + (\gamma -1 )L_0(N) -L_1(N) + E^{(1)}(N)\, .
\eeq

\noindent
CASE $q=2$. Using the representation (\ref{G2eval}),
$$
G^{(2)}_{m,n} = \frac{1}{n} \left(K_2 + K_1\, \log\frac{n}{m} + \frac{1}{4} \log^2\frac{n}{m} \right) + \frac{1}{m}\, S_2(n/m)\, ,
$$
we similarly get
\beqas
Q_N^{(2)} \ceq \sumnl'_{m,n} \left(K_2 + K_1\, \log n - K_1\, \log m + \frac{1}{4} \log^2 n -\frac{1}{2}\log m \log n + \frac{1}{4} \log^2 m\right) a_m\, \frac{a_n}{n} \nonumber\\ && \!\!\!  + \ \sumnl'_{m}\, \frac{a_m}{m}\left( \sumnl'_{n}\, a_n\, S_2(n/m)\right)\nonumber\\ \ceq
\left(\sumnl'_{m} \, a_m\right) \left[\, K_2 \sumnl'_{n}\, \frac{a_n}{n} + K_1 \sumnl'_{n}\, \frac{a_n}{n}\log n + \frac{1}{4} \sumnl'_{n}\, \frac{a_n}{n}\log^2 n\, \right] \nonumber\\ && \!\!\!  
+ \  \left(\frac{1}{4} \sumnl'_{m} \, a_m \log^2 m - K_1 \sumnl'_{m} \, a_m \log m \right)\left(\sumnl'_{n}\frac{a_n}{n}\right) - \frac12 \left(\sumnl'_{m}\, a_m \log m\right)\left(\sumnl'_{n}\frac{a_n}{n} \log n\right)  
\nonumber\\ && \!\!\! 
+ \ \sumnl'_{m}\, \frac{a_m}{m}\, R_2(1/m) + \sumnl'_{m}\, \frac{a_m}{m} \left(\sumnl'_{n}\, a_n\, S_2(n/m) -  R_2(1/m)\right) 
\nonumber\\ \ceq
M_0(N) \left(K_2 L_0(N)  + K_1 L_1(N) + \frac{1}{4} L_2(N)\right) + L_0(N) \left(\frac{1}{4} M_2(N) - K_1 M_1(N) \right)- \frac12 L_1(N) M_1(N)
\nonumber\\ && \!\!\! 
+ \ \sumnl'_{m}\, \frac{a_m}{m}\, R_2(1/m) + E^{(2)}(N)\, .
\eeqas
Since
$$
R_2(x) = \frac{1}{x} \left(\frac12 \log 2\pi + 1 + \frac12 \log x\right) + (2-\gamma)\, \log x - \frac12 \log^2 x + 2\gamma + \gamma_1 - 3
$$ 
for $x\le 1$ and $\frac12 \log 2\pi + 1 = K_1 +\frac{\gamma}{2}$, the fourth term equals
\beqas
\sumnl'_{m} \frac{a_m}{m}\,R_2(1/m) \ceq \sumnl'_{m} \frac{a_m}{m}\left[\, m \left\{\frac12 \log 2\pi + 1 - \frac12 \log m\right\} - (2-\gamma)\, \log m - \frac12 \log^2 m + 2\gamma + \gamma_1 - 3\, \right]\nonumber\\ \ceq
M_0(N) \left(K_1 +\frac{\gamma}{2}\right) - \frac12 M_1(N) - (2-\gamma) L_1(N) -\frac12 L_2(N)  + (2\gamma + \gamma_1 - 3) L_0(N) \, .
\eeqas
Thus, putting things together we obtain
\beqa
Q_N^{(2)} \ceq M_0(N) \left(K_2 L_0(N)  + K_1 \{L_1(N) +1\} + \frac{1}{4}\, \{2\gamma + L_2(N)\} \right) \nonumber\\ && \!\!\! 
- \frac12\,  M_1(N)\,  \bigg( 2K_1\, L_0(N) + L_1(N) + 1\bigg) + \frac{1}{4} M_2(N)L_0(N) 
\nonumber\\ && \!\!\! 
+ \big(2\gamma + \gamma_1 - 3\big)\, L_0(N) + \frac12\, \bigg(2 \gamma L_1(N) - L_2(N)\bigg) - 2 L_1(N) + E^{(2)}(N)\, .\label{Q2eval}
\eeqa

The rationale behind the arrangement of the single terms in these decompositions is twofold. Suppose the coeffients $a_n$ are chosen such that the limits (\ref{L0lim}), (\ref{L1lim}), (\ref{L2lim}) hold also if $L_{\mu,j}$ is replaced by $L_{a,j}$. Then the expressions 
$$
(\gamma -1 )L_{a,0}(N) -L_{a,1}(N)\quad \mbox{and} \quad \big(2\gamma + \gamma_1 - 3\big)\, L_{a,0}(N) + \frac12\, \big(2 \gamma L_{a,1}(N) - L_{a,2}(N)\big) - 2 L_{a,1}(N),
$$
converge to 1 and to 2, respectively, that is, to the desired limit values of the quadratic forms $Q_N^{(q)}$.

To state our second motive for the above decompositions it is convenient to introduce some notation. Here again we temporarily omit the reference to $a$ and simply write $L_j=L_{a,j},\ M_j=M_{a,j}$ etc.
\beqas
P_0^{(1)}(N) \ceq M_0(N)\,  \big[ L_0(N) + \frac12\, \{L_1(N) + 1\}\big],\qquad P_1^{(1)}(N) = M_1(N)L_0(N),\\
P_0^{(2)}(N) \ceq M_0(N)\,  \big[K_2 L_0(N)  + K_1 \{L_1(N) +1\} + \frac{1}{4}\, \{2\gamma + L_2(N)\} \big],\\
P_1^{(2)}(N) \ceq M_1(N)\,  \big[2K_1\, L_0(N) + L_1(N) + 1\big],\qquad P_2^{(2)}(N) 
= M_2(N)L_0(N).
\eeqas
The point is that in these products the factors of the terms $M_j(N)$ are compositions of terms $L_j(N)$ that appear in centered form only: $L_0(N)\equiv \overline L_0(N),\, L_1(N) + 1\equiv \overline L_1(N),\, L_2(N)+2\gamma\equiv \overline L_2(N)$ converge to zero each. This will not suffice to make the products $P_j^{(q)}(N)$ asymptotically negligible. 
However, a new prospect may emerge from the empirical observation that in case of the M\"obius coefficients, $a_n=\mu_n$, the sequences $\overline L_j(N),\, M_j(N),\ N \ge 1$ exhibit an extremely high correlation across different $j$. In fact, the suitably rescaled expressions 
$\widetilde L_j(N) = \overline L_j(N)/\log^j N,\ \widetilde M_j(N) = M_j(N)/\log^j N$ are visibly almost indistinguishable for different $j$; see Appendix, Fig.~2. One may therefore expect that their differences are much reduced. Indeed, Fig.~3 shows a massive size reduction when passing to the latter.

\smallskip
An idea how this observation might be put to use is to ask for coefficients $a_n$ which are such that the {\em raw} quantities $\overline L_{a,j}(N),\, M_{a,j}(N)$ built with the coefficients $a_n$ mimic the behavior of the {\em differences}
\beq\label{Deldef0}
\Delta \overline L_{\mu,j}(N) = \overline L_{\mu,j}(N) - \overline L_{\mu,j+1}(N)/\log N\, ,\qquad \Delta M_{\mu,j}(N) = M_{\mu,j}(N)-M_{\mu,j+1}(N)/\log N
\eeq
built with the M\"obius coefficients $\mu_n$. The next subsection answers this question in the affirmative.

\subsection{Choice of the coeffients}\label{chco}

We begin by introducing some notation. Let $\kappa(s) = 1/\zeta(s)$, and put
$$
\eta_j = (-1)^j \kappa^{(j)}(2) = \sumnl_{n\ge 1} \frac{\mu_n}{n^2}\,\log^j n\, , \ j = 0,1,\ldots .
$$
Recall that for given coefficients $a_n = a_{n,N}$ we denote
\beqas
M_{a,j}(N) \ceq \snn a_{n,N}\, \log^j n\, , \quad L_{a,j}(N) = \snn \frac{a_{n,N}}{n}\, \log^j n\, \quad \left(\snn \equiv \sumnl_{n \le N}\right) \qquad \mbox{and}\\
\overline L_{a,j}(N) \ceq L_{a,j}(N) + \ell_j\quad (j = 0,1,\ldots), \quad \mbox{where}\quad \ell_0 =0,\, \ell_1 =1,\,\ell_2 =2\gamma,\,\ell_3 =6(\gamma^2 + \gamma_1) .
\eeqas
We construct the searched-for coefficients $a_{n,N}$ as a modification $\nu_{n,N}$ of the Levinson-Selberg coefficients $\lambda_{n,N} = \mu_n \left(1 - \log n/ \log N\right)$, to wit,
\beq\label{nudef}
\nu_{n,N} = \mu_n \left(1 - \frac{1}{\log N}\, \left[\log n + a\frac{\log n}{n} + \frac{b}{n}\right]\right) - \frac{c}{n\log N}
\eeq
with certain constants $a,b,c$.

\begin{prop}\label{nucoef}
Suppose that the constants $a,b,c$ satisfy the system of equations 
\beqas
a \eta_1  + b \eta_0 + c \zeta(2)\ceq 1 \label{eq1}\\
a \eta_2  + b \eta_1 - c \zeta'(2)\ceq 2\gamma\label{eq2}\\
a \eta_3  + b \eta_2 + c \zeta''(2)\ceq 6(\gamma^2 + \gamma_1)\, .\label{eq3}
\eeqas
Then for $j = 0 ,1,2$
\beqa
\overline L_{\nu,j}(N) \ceq \Delta \overline L_{\mu,j}(N) + O(N^{-1} \log^j N)\, ,\label{Lappr}\\
M_{\nu,j}(N) \ceq \Delta M_{\mu,j}(N) + O(\log^j N)\, .\label{Mappr}
\eeqa
\end{prop}
{\bf {\em Proof.} }  It suffices to consider one of the three cases $j=0,1,2$. E.g., for $j=1$ we get
\beqas
&& \!\!\!\! \overline L_{\nu,1}(N)  = L_{\nu,1}(N) + 1 = \snn \frac{\nu_{n,N}}{n}\log n\ +\  1 \nonumber\\ \ceq L_{\mu,1}(N) + 1
 - \frac{1}{\log N}\, \left[\snn  \frac{\mu_n \log^2 n}{n} + a \snn  \frac{\mu_n \log^2 n}{n^2} + b \snn  \frac{\mu_n\log n}{n^2} + c\snn  \frac{\log n}{n^2}\right]\nonumber\\ \ceq
\overline L_{\mu,1}(N)  - \frac{1}{\log N}\, \left[L_{\mu,2}(N) + a \eta_2 + b \eta_1 - c\zeta'(2) + O(N^{-1}\log^2 N) \right]\, ,
\eeqas
using the straightforward estimate $\sumnl_{n > N} n^{-2} \log^j n = O(N^{-1} \log^j N)$ and the definition of the $\eta_j$s.  By the second equation of the system the last line reduces to $\overline L_{\mu,1}(N)  - \overline L_{\mu,2}(N)/\log N + O(N^{-1}\log N)$, which completes the proof of (\ref{Lappr}) for $j=1$. \\
The $M$-terms can be dealt with analogously. If $j=2$, say, then
\beqas
&& \!\!\!\! M_{\nu,2}(N) = \snn  \nu_{n,N}\log^2 n \nonumber\\ \ceq 
M_{\mu,2}(N)
 - \frac{1}{\log N}\, \left[\snn  \mu_n \log^3 n + a \snn  \frac{\mu_n \log^3 n}{n} + b \snn  \frac{\mu_n\log^2 n}{n} + c\snn  \frac{\log^2 n}{n} \right]\nonumber\\ \ceq
M_{\mu,2}(N)  - \frac{1}{\log N}\, \left[M_{\mu,3}(N) + a L_{\mu,3}(N) + b L_{\mu,2}(N) + O(\log^3 N) \right]\nonumber\\ \ceq
M_{\mu,2}(N)  - \frac{M_{\mu,3}(N)}{\log N} + O(\log^2 N)\, ,
\eeqas
the latter since $L_{\mu,3}(N) = o(\log N)$; cf.~end of proof of Proposition \ref{mxlim}.
\done

\smallskip\noindent
{\bf Remarks.} 1. A Taylor expansion of $\frac{1}{\zeta(s)}$ at $s=1$ suggests that $\lim_{N\to\infty} L_{\mu,3}(N) = -6(\gamma^2\! +\! \gamma_1)$, yet a corresponding reference is not known to this author.  
In fact, our argument goes without this convergence. \\
2. The limits $\lim_{N\to\infty} L_{\nu,j}(N),\, j=0,1,2$ clearly are identical to those of $L_{\mu,j}(N)$. Likewise, Proposition \ref{mxlim} continues to hold if the coefficients $\lambda_{n,N}$ are replaced by $\nu_{n,N}$. \\
3. Explicitly known in the above system of equations are only the entries $\zeta(2)=\pi^2/6$ and $\eta_0 = 1/\zeta(2)$. Numerical evaluation yields the approximative solution
$$%
a \doteq -2.116586,\quad b \doteq -0.407487,\quad c \doteq 0.312679\, .
$$
If  $q=1$, only the cases $j=0,1$ are of interest, and one may omit one of the additional entries in (\ref{nudef}); the last one, say. The values of the solution $a,b$ of the reduced system then differ, of course.  \\
4. The centering of the $L$-terms as achieved by the additional terms in (\ref{nudef}) is essential for the remarkable match stated in (\ref{Lappr}). The Levinson-Selberg coefficients $\lambda_{n,N}$ would not do. This is different with the $M$-terms, which are relatively little affected by such a modifcation.

\subsection{Estimating the $L$-and $M$-terms: Two notes}\label{LMest}

The approach sketched above hinges on good estimates for the $L$-and $M$-terms. Until present these seem to be very far away. Let us nonetheless give a brief account of two methods that have been applied in this context. Throughout the sequel the coefficients $a_n$ of the $L$-and $M$-terms are supposed to be given by the M\"obius function, $a_n = \mu(n)$. In return, our target objects will be the differences $\Delta \overline L_\mu,\ \Delta M_\mu$ which by Proposition \ref{nucoef} may serve as approximations to $\overline L_\nu\, ,M_\nu$. 

\medskip\noindent
{\bf 1. MacLeod type identities and M\"obius inversion.} Our presentation follows M.~Balazard's ``Remarques élementaires'' \cite{Ba}. His starting point is the following basic observation.

\medskip\noindent
\cite[Eq.~(6)]{Ba} {\em Let $f(n)$ be an arithmetic function, $F(x) = \snx f(n)$, and suppose that $\phi(x)$ is absolutely continuous for $x\ge 1$. Then}
\beqa
\Phi(x)\equiv  
\snx f(n) \phi(x/n) 
 \ceq \int_1^x F(x/t)\, \phi'(t)\, dt + F(x) \phi(1) \quad  (x \ge 1).\label{f1}
\eeqa

The proposition allows to derive bounds for the function $\Phi(x)$  in terms of bounds on the function $F(x)$. Producing related identities of interest requires some ingenuity; e.g. \cite{Ba,D,R23}. Balazard works with identities involving the fractional part function due to MacLeod \cite{ML}, and  with M\"obius inversion. Here we take the functions $R_1(x),\, R_2(x)$ defined in Theorem \ref{mainthm} as our starting points. The two propositions below are closely related to, respectively, Balazard's \cite{Ba} Proposition 10 and to Ramaré's \cite{R15} Theorem 1.7 and Lemma 5.1; see also F. Daval \cite[Chapitre 6]{D}.

\begin{prop}  \label{pid1}
Let $\phi_1(x) = xR_1(x)$, and set $\Phi_1(x) = \snx \mu(n) \phi_1(x/n)$. Then 
\beq\label{lmrel1}
(1-\gamma) L_0(x) - \frac12\, \frac{M_0(x)}{x} = \Delta \overline L_0(x)\, \log x  - \frac{\Phi_1(x)}{x}  + \frac{1}{x} \qquad (x \ge 1).
\eeq
\end{prop}
{\bf {\em Proof.} } Let us begin by recalling the formula (\ref{R1def}),
$$
R_1(x) = \log x + \gamma - H(x) - \frac{\{x\} - 1/2}{x} = \log x + \gamma -1 - H(x) + \frac{\lfloor x \rfloor + 1/2}{x}\, .
$$
Let $w(x) = (x-1)(x\ge 1)$, where $(x \in A)$ denotes the indicator function of  $A$, and
put $W(x) = \snx w(x/n)$. Then
\beqas
W(x) = \snx \left(\frac{x}{n}-1\right)(x/n \ge  1) = x H(x) - \lfloor x \rfloor =
 x\log x + (\gamma - 1)x + 1/2 - xR_1(x).
\eeqas
M\"obius inversion gives, for $x \ge 1$,
\beqas
x-1 \ceq \snx \mu(n)\, W(x/n) \\ \ceq
\snx \mu(n) \frac{x}{n}\log \frac{x}{n} + (\gamma-1) \snx \mu(n) \frac{x}{n} + \frac12\snx \mu(n) - \snx \mu(n) \phi_1(x/n)\\ \ceq
L_0(x)\, x \log x - x \left(L_1(x)+1\right) + x + (\gamma - 1)x L_0(x) + \frac12\, M_0(x) - \Phi_1(x) \, .
\eeqas
Subtracting $x$, then dividing by $x$ we get 
\beqas
-\frac{1}{x} \ceq \Delta \overline L_0(x)\, \log x  + (\gamma - 1) L_0(x) + \frac12\, \frac{M_0(x)}{x} - \frac{\Phi_1(x)}{x} \, ,
\eeqas
which is  (\ref{lmrel1}). \done

\medskip
\begin{prop} \label{pid2}
Let $\phi_2(x) = xR_2(x)$, and set $\Phi_2(x) = \snx \mu(n) \phi_2(x/n)$.\\ Moreover, let $C_1 = 2\gamma + \gamma_1 - 3,\ C_2 = \frac12 \log 2\pi + 1$. Then
\beqa\label{lmrel2}
C_1 L_0(x) + C_2\, \frac{M_0(x)}{x} \ceq \frac{\log x}{2} \left(\Delta \overline L_0(x)\, \log x - \Delta \overline L_1(x) -  \frac{\Delta M_0(x)}{x} \right) -(2-\gamma) \Delta\overline L_0(x) \log x \nonumber\\ && 
+ \frac{\Phi_2(x)}{x} + O\left(\frac{\log x}{x}\right) \qquad (x \ge 1).
\eeqa
\end{prop}
{\bf {\em Proof.} }
Let $W(x) = \snx w(x/n)$ where $w(x) = (1+x)\log x + 2 -2x$. 
The representation (\ref{R2def}) of $R_2(x)$ can then be rewritten as
$$
W(x) = \frac{x}{2} \log^2 x - (2-\gamma)\,x \log x -(2\gamma + \gamma_1 - 3)x - \frac12 \log x - \frac12 \log 2\pi - 1 + xR_2(x)\, .
$$
M\"obius inversion gives, for $x \ge 1$,
\beqas
w(x) \ceq \snx \mu(n) W(x/n) = x\snx \frac{\mu(n)}{n} \left[ \frac12 \log^2\frac{x}{n}  - (2-\gamma)\log \frac{x}{n}  - C_1\right] \\
&& \ - \frac12\snx \mu(n)\log \frac{x}{n} - C_2 \snx \mu(n) + \snx \mu(n) \phi_2(x/n) \\ \ceq
\frac{x}{2} \left\{L_0(x)\log^2 x - 2 L_1(x)\log x + L_2(x)\right\} - (2-\gamma) x \left\{L_0(x)\log x - L_1(x)\right\} \\
&&  - C_1 x L_0(x) - \frac12 M_0(x)\log x + \frac12 M_1(x)  - C_2M_0(x) + \Phi_2(x)
\\ \ceq
\frac{x}{2} \left\{L_0(x)\log^2 x - \overline L_1(x) \log x - \overline L_1(x) \log x + \overline L_2(x) + 2\log x - 2\gamma\right\} 
\\ && 
-\ (2-\gamma) x \left\{L_0(x)\log x - \overline L_1(x) + 1\right\} - \frac12 \Delta M_0(x)\log x  - C_1 x L_0(x) - C_2M_0(x)  + \Phi_2(x)
\eeqas
After division by $x$ this becomes
\beqas
\log x - 2 + O\left(\frac{\log x}{x}\right) \ceq \frac12 \left\{\Delta \overline L_0(x) \log^2 x  - \Delta \overline L_1(x) \log x\right\} + \log x - \gamma - (2-\gamma) \Delta \overline L_0(x)\log x
\\ && \!\!\!
-\ 2 + \gamma  - \frac12 \frac{\Delta M_0(x)}{x}\log x - C_1 L_0(x) - C_2\, \frac{M_0(x)}{x} + \frac{\Phi_2(x)}{x}\, .
\eeqas
With the appropriate cancellations we get
\beqas
O\left(\frac{\log x}{x}\right) \ceq 
\frac{\log x}{2} \left(\Delta \overline L_0(x)\, \log x - \Delta \overline L_1(x) -  \frac{\Delta M_0(x)}{x} \right) -(2-\gamma) \Delta\overline L_0(x) \log x 
\\ && \!\!\!
- C_1 L_0(x) - C_2\, \frac{M_0(x)}{x}  + \frac{\Phi_2(x)}{x}\, ,
\eeqas
and the claim follows. \done

\newpage
\medskip
Other identities follow by linear combination. For example, if one subtracts equation (\ref{lmrel2}) from equation (\ref{lmrel1}) multiplied by $(C1+C2)/(1/2-\gamma)$, the left-hand side becomes a multiple of $L_0(x) - M_0(x)/x$, which quantity is a  central study object in \cite{Ba}. By similar linear combinations one can delete the term $L_0$ at the left, thus isolating the term $M_0$; and vice versa. (Recall that $M_0(x)$ is the Mertens function.) 
Estimates for the terms  $\Phi_q(x)$ can be obtained as described in \cite{Ba}. Proposition \ref{RS} combined with, respectively, (\ref{R1rep}) and (\ref{R2rep2}), (\ref{Vrep}) entails 
$$
|\phi_1'(x)| = O(1/x), \qquad |\phi_2'(x)| = O(1/x^2)\qquad (x \to\infty).
$$
An application of (\ref{f1}) in turn yields the bounds 
$$
|\Phi_1(x)| \le A_1 \left(\int_1^x \frac{|M_0(t)|}{t}\, dt +  |R_1(1)\, M_0(x)|\right) \, , \quad |\Phi_2(x)| \le A_2 \left(\frac{1}{x} \int_1^x |M_0(t)|\, dt + |R_2(1)\, M_0(x)|\right) 
$$
where $A_1,A_2$ are finite constants, and $R_1(1)=\gamma-1/2,\, R_2(1)=C_1+C_2\doteq .000554$.
Improved estimates can be derived by skillful variations of the identity (\ref{f1}), cf.~\cite{Ba,R15,D,D1}. The motive for the present choice was to feature the link leading from the functions $R_1\, , R_2$ to the difference terms $\Delta \overline L_j\, , \Delta M_j$. 

\medskip
The literature on the matter touched on here is vast. We may once more point to the work of  Balazard \cite{Ba}, Ramaré \cite{R15}, and Daval \cite{D} for references and valuable accounts of the history of the subject.
Recently, Ramaré and Zuniga-Alterman \cite{R23} unified the matter by substituting the factors $n^{-1}$ and $n^0\equiv 1$ figuring in the definition of the $L$- and $M$-terms, respectively, by the general factor $n^{-s}$, where $s$ is a complex parameter. See also the earlier paper \cite{R13} as well as \cite{D1}.

\bigskip\noindent
{\bf 2. Perron's formula.} 
Let us briefly comment on this approach on the basis of the example 
\beqa
\Delta \overline L_1(N) \log N \ceq 
\snn \frac{\mu_n \log n}{n} \log \frac{N}{n} + \log N -2\gamma \, .\label{L1xpd}
\eeqa
We essentially follow Bettin et al. \cite{BCF}. The Dirichlet series associated with the coefficients $n^{-1}\mu_n \log n$ is $\zeta'(s+1)/\zeta(s+1)^2$. 
The appropriate Perron's formula thus is
\beq\label{perron1}
2\pi i\, \snn \frac{\mu_n \log n}{n} \log \frac{N}{n} = \int_c \frac{N^{w-1}\, \zeta'(w)}{(w-1)^2\, \zeta(w)^2}\, dw = \int_c \frac{N^{w-1}}{(w-1)^2}\left[\log N - \frac{2}{w-1}\right] \frac{1}{\zeta(w)}\, dw
\eeq
where again $\int_c$ denotes integration upwards the vertical line from $c - i\infty$ to $c + i \infty$, and $c>1$. 
Let $G_N(w)$ denote the integrand at the right-hand side of (\ref{perron1}). ``Now we move the path of integration to $\Re w = -2M -1$ for some large integer $M$'' \cite{BCF} 
to obtain, formally at least, 
\beq\label{Geval}
\frac{1}{2\pi i} \int_c G_N(w)\, dw = \R\ G_N(w) \big|_{w=1} + \sumnl_\rho \R\ G_N(w) \big|_{w=\rho} + \sumnl_{n \ge 1} \R\ G_N(w) \big|_{w=-2n} 
\eeq
where the second and the third sum of residues extend over the critical zeros $\rho$, and the trivial zeros $-2,\,-4,\ldots$ of $\zeta$, respectively. The residue $\R\ G_N(w) \big|_{w=1} = -\log N + 2\gamma$ cancels the offset in (\ref{L1xpd}). For $v$ a simple zero of the zeta function, the residue 
$$
\R\ G_N(w) \big|_{w=v} = \frac{N^{v-1}}{\zeta'(v)\, (v-1)^2}\left[\log N - \frac{2}{v-1}\right] \qquad (v \neq 1).
$$
Using this formula along with a well-known expression for $1/\zeta'(-2n)$, the last sum in (\ref{Geval}) can be shown to be of the order $O(N^{-3}\log N)$ \cite{BCF}, so if we assume that all $\zeta$-zeros are simple we obtain
\beqas
\Delta \overline L_1(N) \log N \ceq \sumnl_{\rho} \frac{N^{\rho-1}}{\zeta'(\rho)\, (1-\rho)^2}\left[\log N + \frac{2}{1-\rho}\right] + O\!\left(\frac{\log N}{N}\right) .
\eeqas
This is essentially  $N^{-1} \log N$ times the sum of residues $\sumnl_{\rho} \frac{N^{\rho}}{\zeta'(\rho)\, \rho^2}$ dealt with in \cite{BCF}.

\newpage

\small
\bibliographystyle{amsalpha}

\smallskip\noindent
E-mail address: wernehm@web.de

\normalsize

\newpage

\section*{Appendix}\label{plots}

In this appendix we present some plots for illustration. 
Figure 1 depicts the functions $R_q,\, S_q,\ q=1,2$. For $x \ge 3/4$, the series $S_q(x) = \sum_n R_q(nx)$ is largely dominated by its first term, $R_q(x)$. At a finer scale $S_1(x)$ exhibits myriads of little spikes occurring at rational $x$. The left-hand boundary behavior of $S_q(x)$ follows from Proposition \ref{RS} on making use of the reciprocity relations  (Corollary \ref{reci}): 
as $x \to 0$
$$
xS_1(x) = \frac12\, |\log  x| - K+ o(1),\qquad xS_2(x) = -\frac{1}{4}\log^2 x + K_1\, |\log x| - K_2+ o(1).
$$

\vspace*{-2mm}
\begin{minipage}{.44\linewidth}
\includegraphics[scale=.28,angle=0,clip]{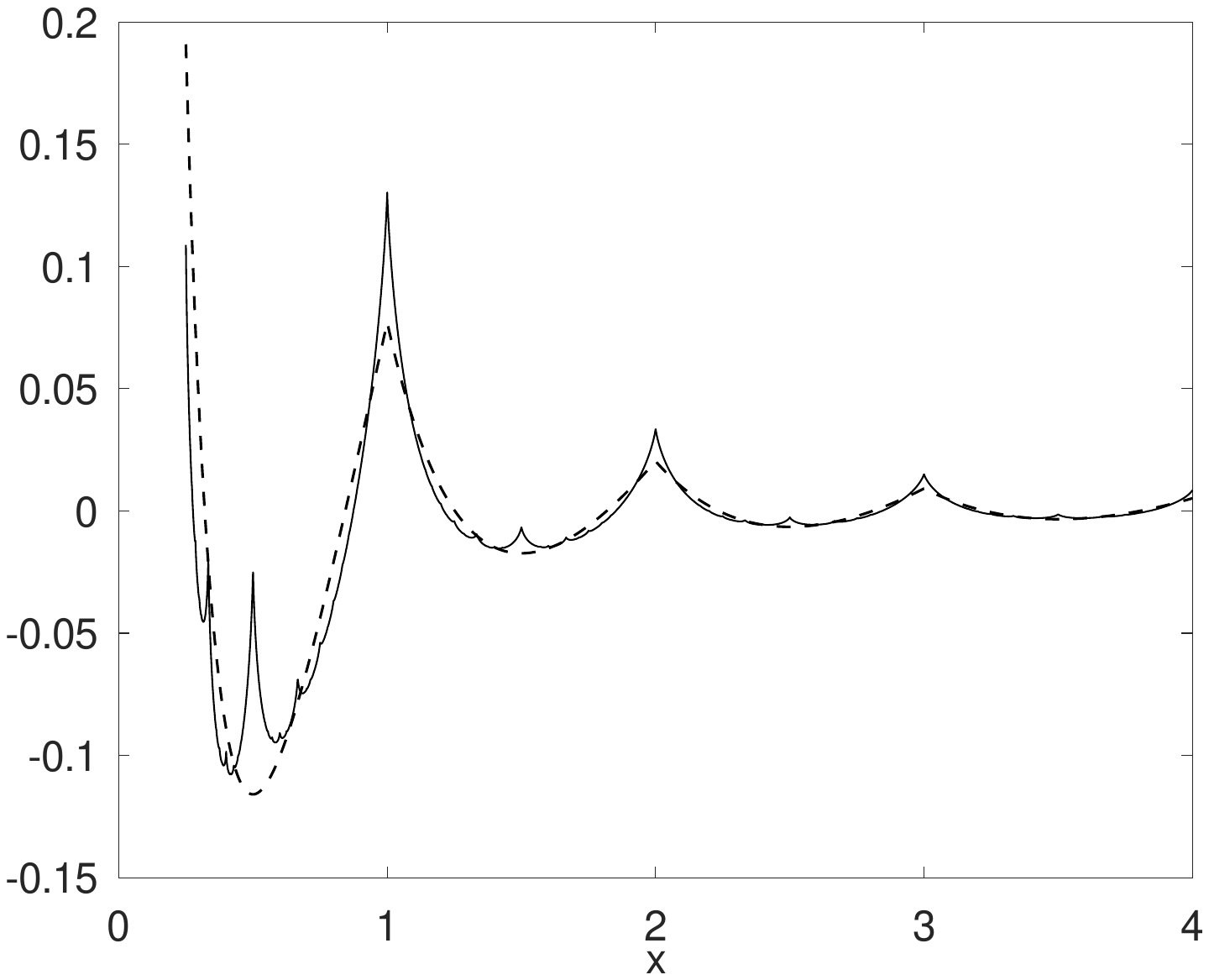}
\end{minipage} 
\hspace{8mm} 
\begin{minipage}{.44\linewidth} 
\includegraphics[scale=.28,angle=0,clip]{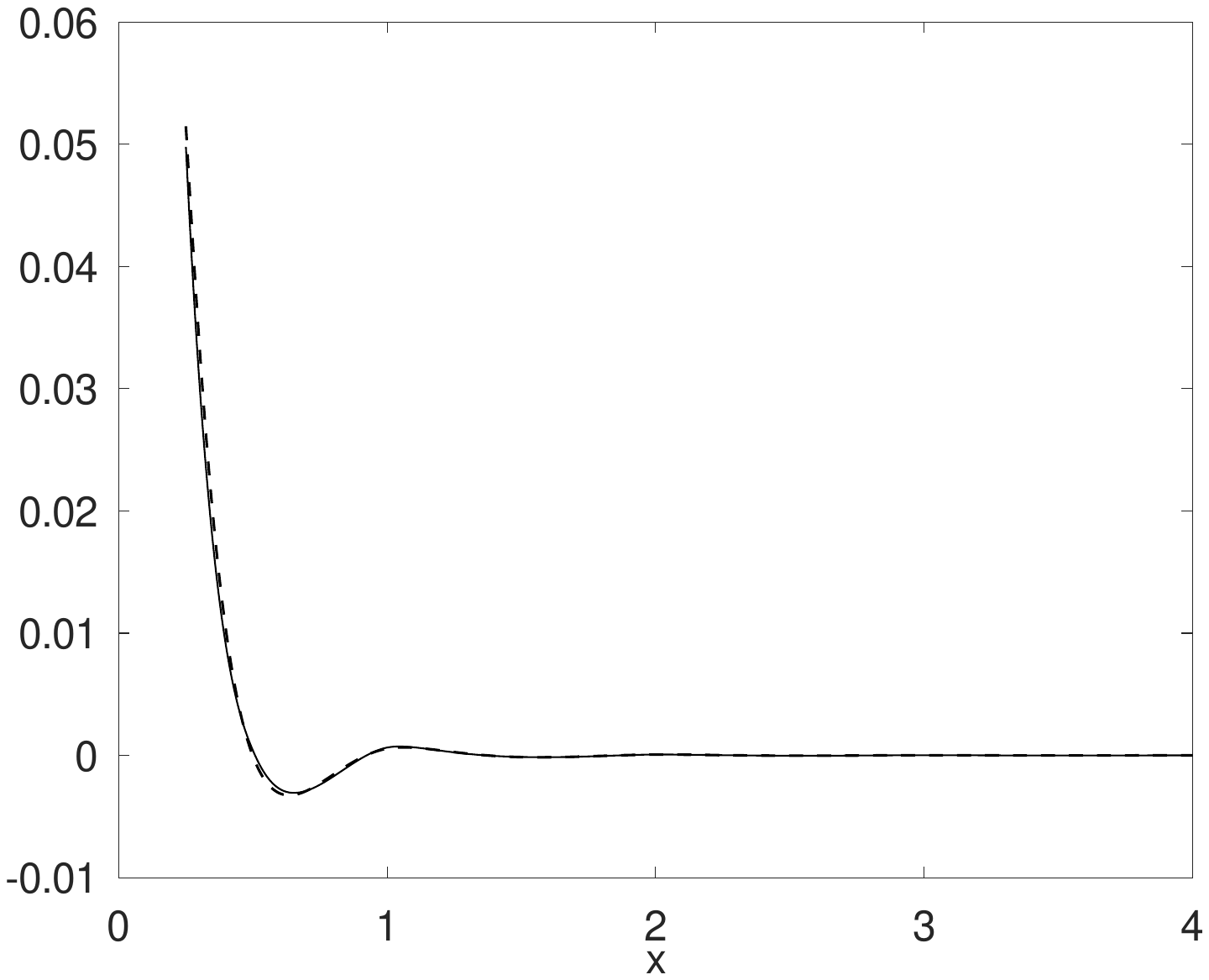}
\end{minipage} \hfill
\newline

\vspace*{-2mm}
Figure 1.  \small{Left: Functions $R_1(x)$ (dashed) and $S_1(x)$ (solid) plotted vs.~$x$ in the range $1/4 \le x \le 4$.\\ Right: Same for $R_2,\, S_2$.}
\normalsize

\medskip
Figure 2 demonstrates the strong correlation phenomenon mentioned at the end of Section \ref{secqdec}. The traces of the rescaled $L$- and $M$-terms $\widetilde L_j(N) = \overline L_j(N)/\log^j N,\ \widetilde M_j(N) = M_j(N)/\log^j N$ dovetail almost perfectly.

\vspace*{4mm}
\begin{minipage}{.44\linewidth}
\includegraphics[scale=.28,angle=0,clip]{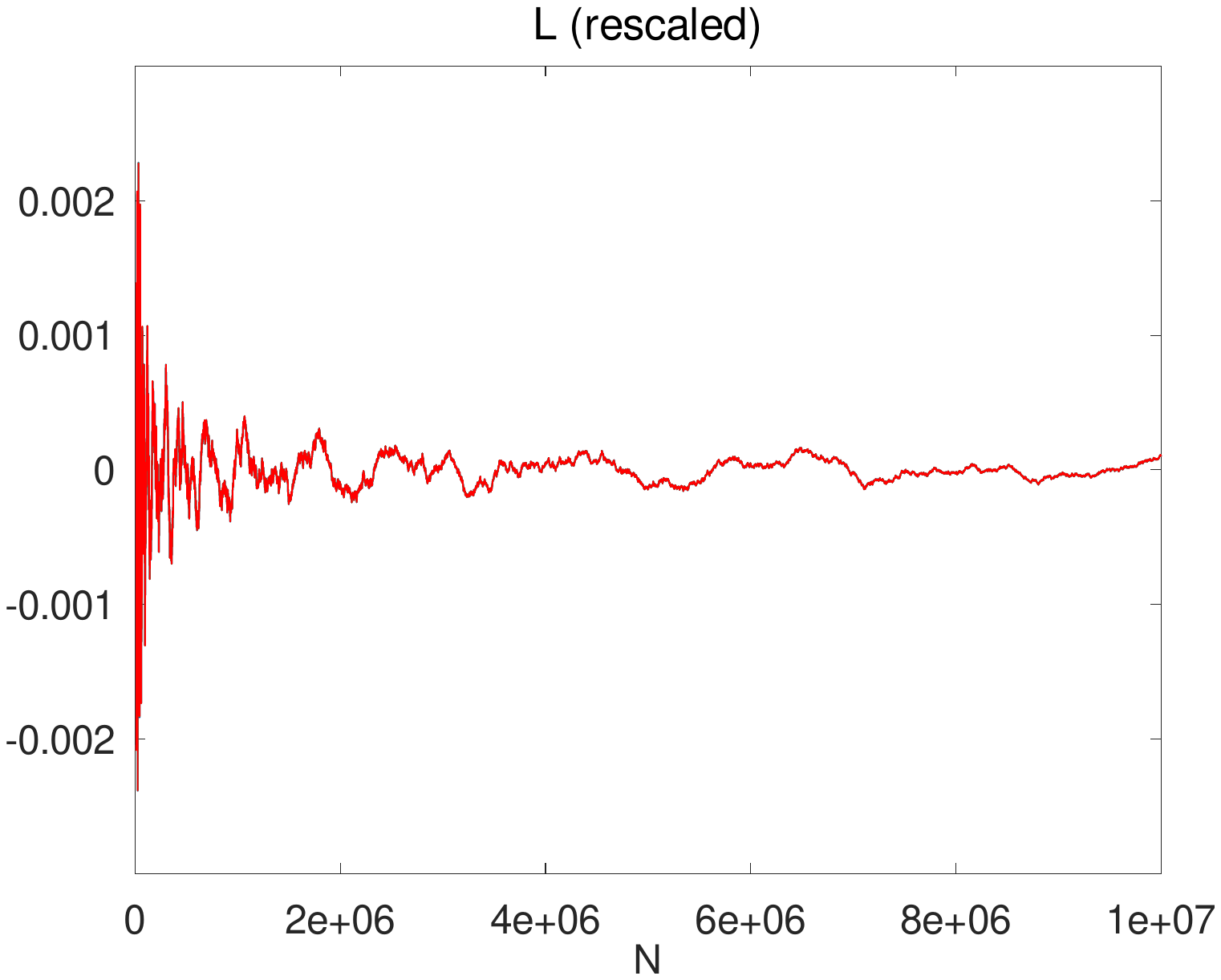}
\end{minipage} 
\hspace{8mm} 
\begin{minipage}{.44\linewidth} 
\includegraphics[scale=.28,angle=0,clip]{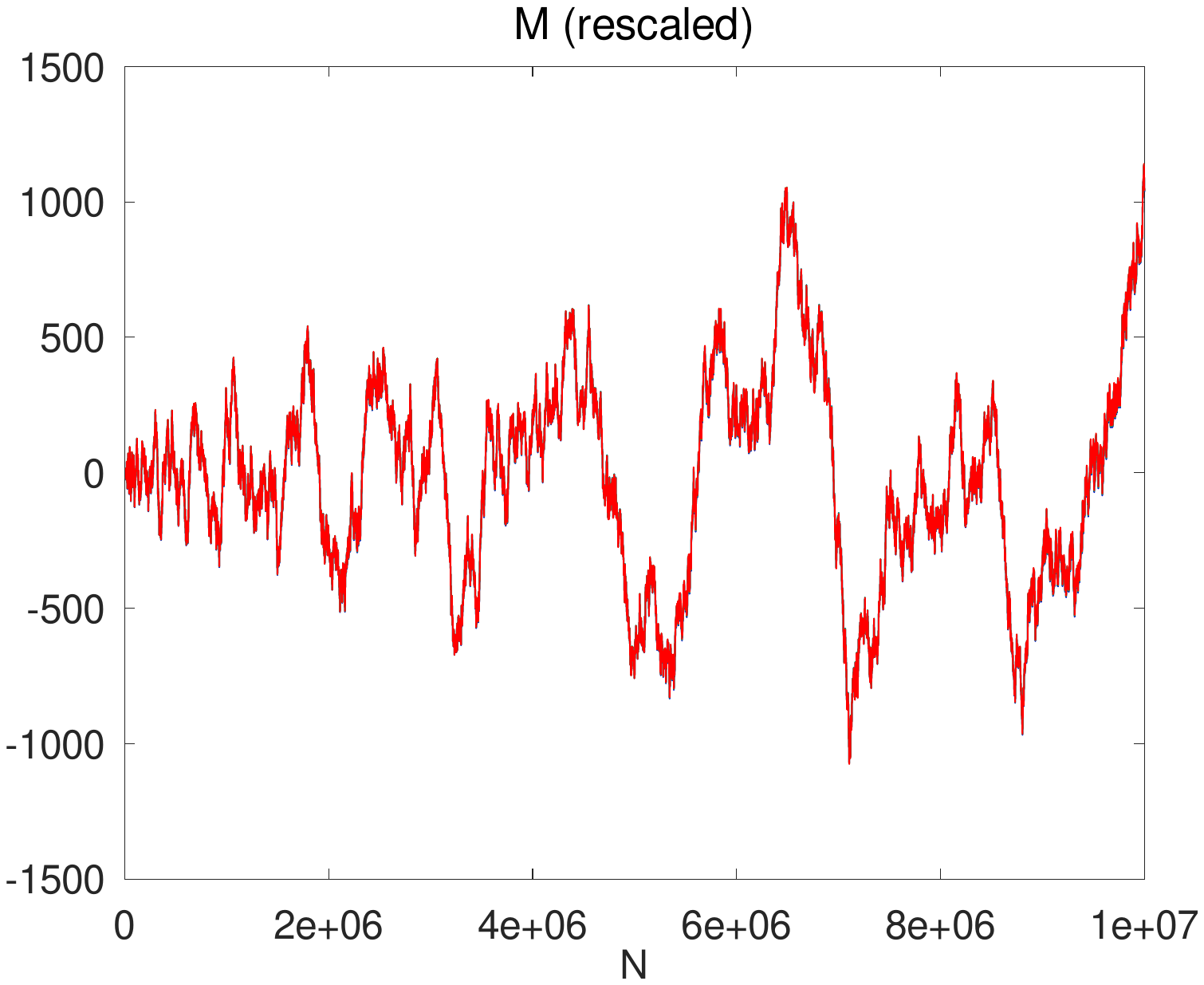}
\end{minipage} \hfill
\newline

\vspace*{-2mm}
Figure 2. \small{
Left: $\widetilde{L}_{\mu,j}(N)$ plotted vs.~$N$ for $j=0$ (blue), $j=1$ (green), $j=2$ (red) in the range $10^4 \le N\le 10^7$ at every 500-th data point. (This does not affect the picture's general features.) Right: Same for $\widetilde{M}_{\mu,j}(N)$. }  
\normalsize

\medskip
Figure 3 highlights the massive reduction of the deflections when passing to the differences of the rescaled terms $\Delta \widetilde{L}_j(N) = \widetilde{L}_j(N) - \widetilde{L}_{j+1}(N)$ etc. Surprisingly, the nearly perfect congruence apparent in Figure 2  even extends to these differences, with one exception. The course of $\Delta \widetilde{M}_{\mu,0}(N)$ in the right panel deviates from the hardly distiguishable $\Delta \widetilde{M}_{\mu,1}(N),\, \Delta \widetilde{M}_{\mu,2}(N)$. The shift vanishes if the former is replaced by $\Delta \widetilde{M}_{\mu,0}(N) + \pi^2/6$ (not shown). 
In accordance with Proposition \ref{nucoef} the courses of the rescaled $L$-terms $\widetilde {L}_{\nu,j}(N)$ built with the coefficients $\nu_{n,N}$ are visibly indistinguishable from the differences $\Delta \widetilde{L}_{\mu,j}(N)$ shown in Figure 3, left panel. This applies also to the corresponding $M$-terms. 

\medskip
Figure 4 addresses the question to what extent the $L$-factors can downsize the $M$-factors in the products figuring in the decompositions (\ref{Q1eval}), (\ref{Q2eval}). The results look encouraging. We specifically consider the ultimately relevant sums of the product terms built with coefficients $\nu_{n,N}$,
$$
P_\nu^{(1)}(N) = P_{\nu,0}^{(1)}(N) + P_{\nu,1}^{(1)}(N), \qquad P_\nu^{(2)}(N) = P_{\nu,0}^{(2)}(N) + P_{\nu,1}^{(2)}(N) + P_{\nu,2}^{(2)}(N)\, .
$$

\begin{minipage}{.44\linewidth}
\includegraphics[scale=.28,angle=0,clip]{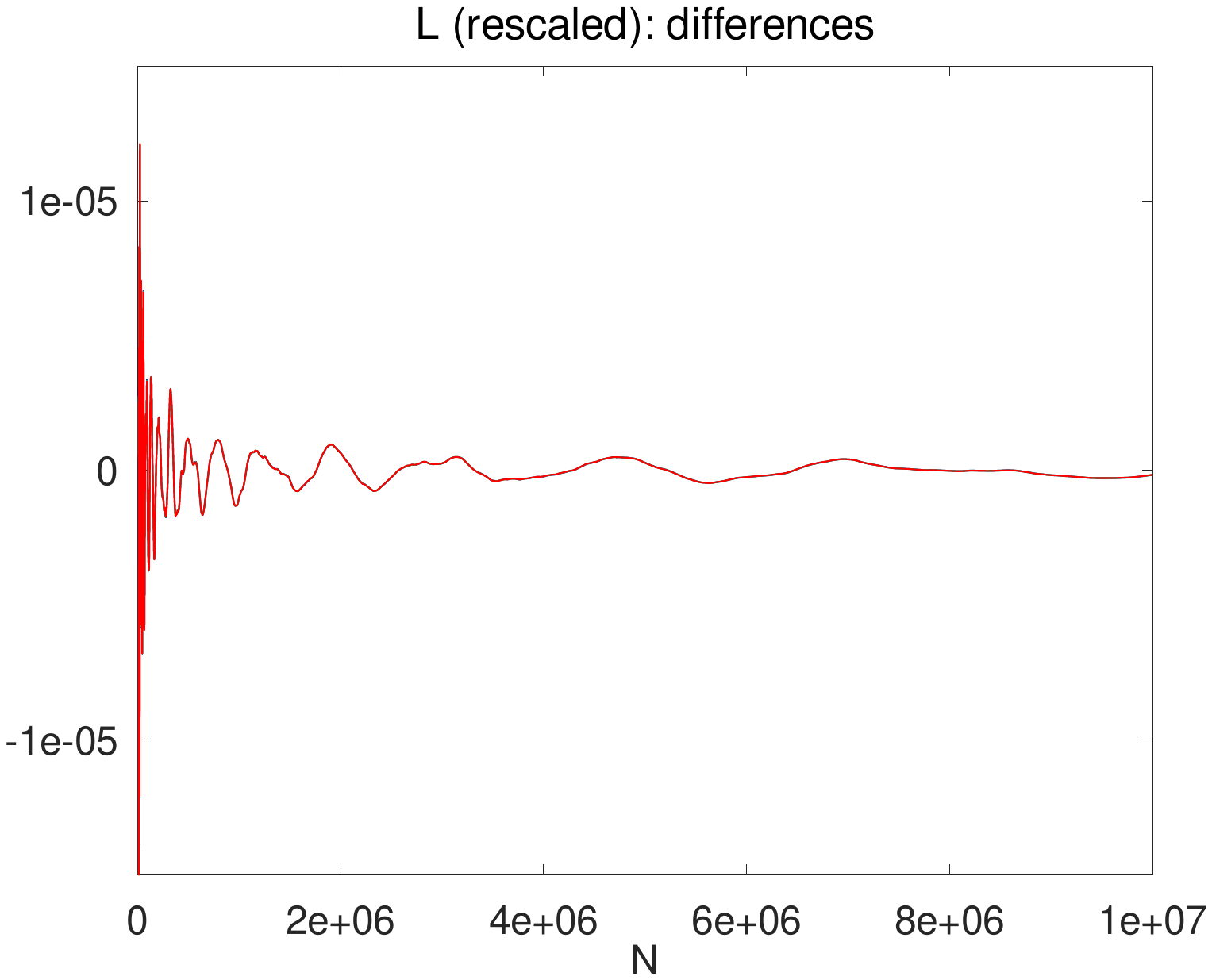}
\end{minipage} 
\hspace{8mm} 
\begin{minipage}{.44\linewidth} 
\includegraphics[scale=.28,angle=0,clip]{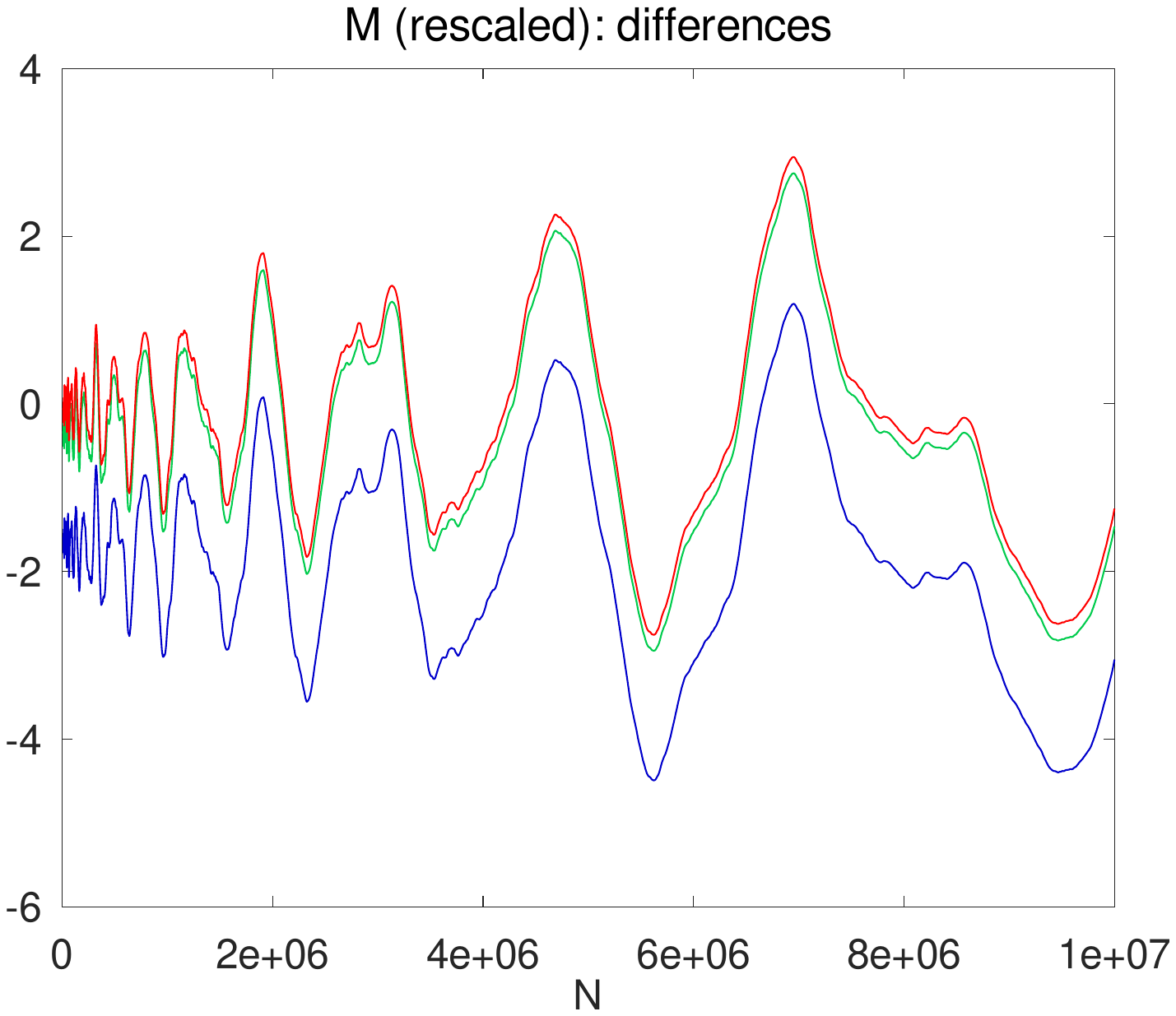}
\end{minipage} \hfill
\newline

\vspace*{-2mm}
Figure 3. \small{
Traces of $\Delta \widetilde{L}_{\mu,j}(N)$ (left) and $\Delta \widetilde{M}_{\mu,j}(N)$ (right); otherwise as in Fig.~2.
}
\normalsize

\bigskip
\begin{minipage}{.44\linewidth}
\includegraphics[scale=.27,angle=0,clip]{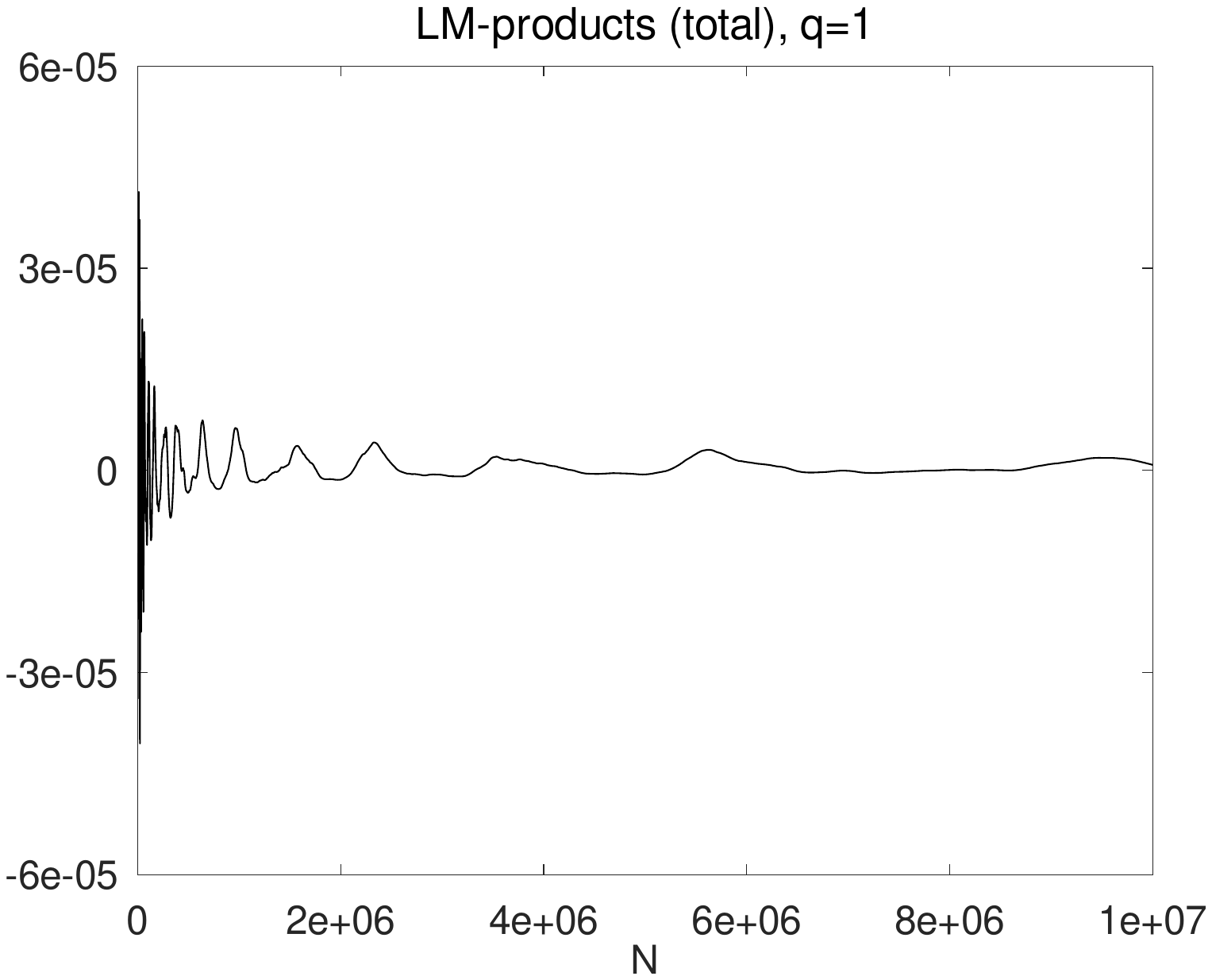}
\end{minipage} 
\hspace{8mm} 
\begin{minipage}{.44\linewidth} 
\includegraphics[scale=.27,angle=0,clip]{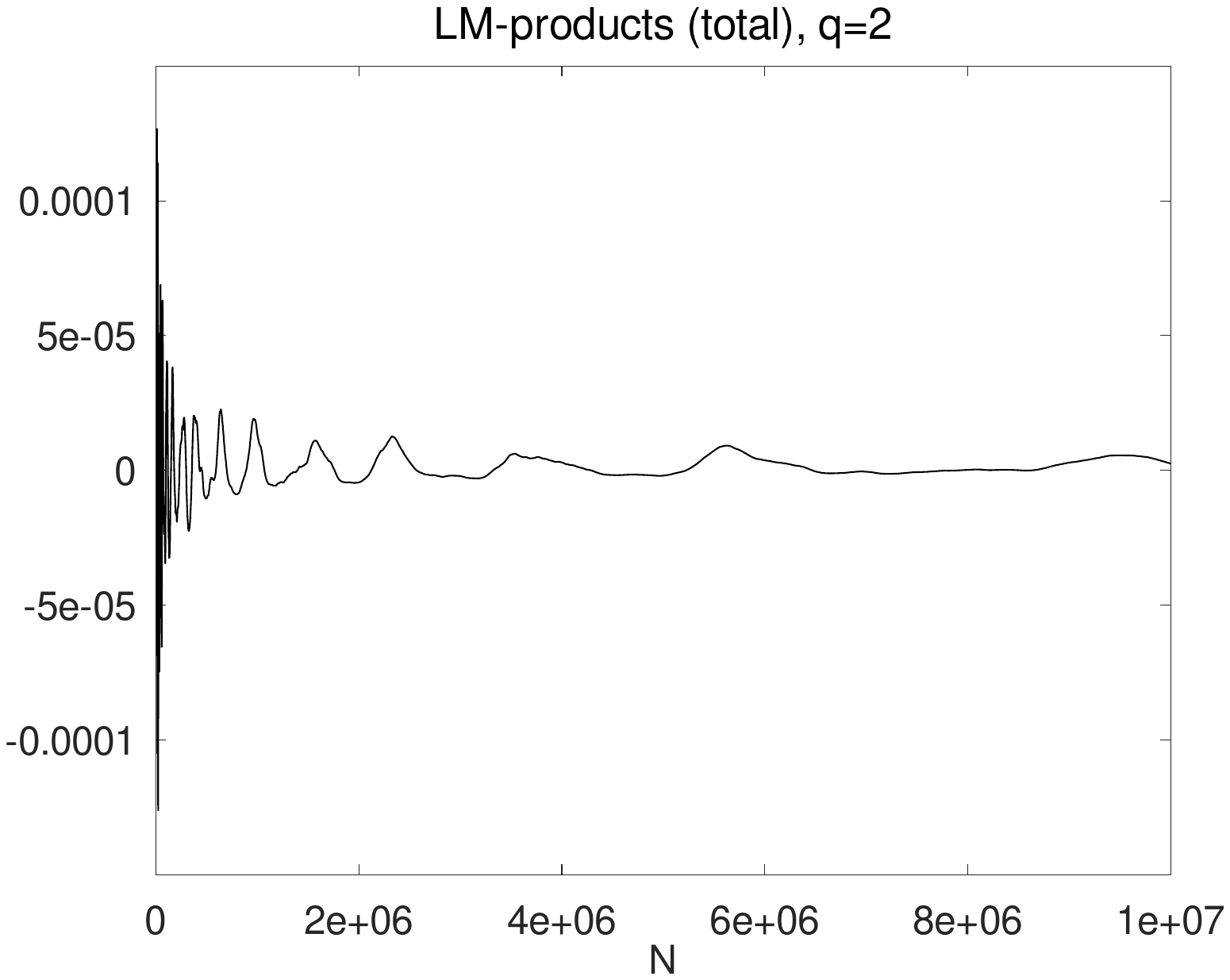}
\end{minipage} \hfill
\newline

\vspace*{-2mm}
Figure 4. \small{
Traces of $P_\nu^{(q)}(N) $ in the range $10^4 \le N\le 10^7$ at every 500-th data point.
}
\normalsize

\medskip
Figure 5 compares scaled distances $d_q(N)$ (cf.~(\ref{cplxerr})) computed with the mollified M\"obius coefficients $\lambda_{n,N}=\mu_n (1-\log n/\log N)$ to those computed with the coefficients $\nu_{n,N}$. 
In the case $q=1$ it is known that under RH the distance (or approximation error) $d_1(N)$ cannot tend to zero faster than $O(1/\sqrt{\log N})$, which is believed to be the correct rate; see \cite{BBLS0,Bu2,LR}. 
The picture at the right suggests a possibly faster decay rate when $q=2$. However, this is weak evidence given the modest range of $N$. Note that the error for the $\nu$-coefficients is persistently smaller than the one for the $\lambda$-coefficients.

\vspace*{2mm}
\begin{minipage}{.44\linewidth}
\includegraphics[scale=.28,angle=0,clip]{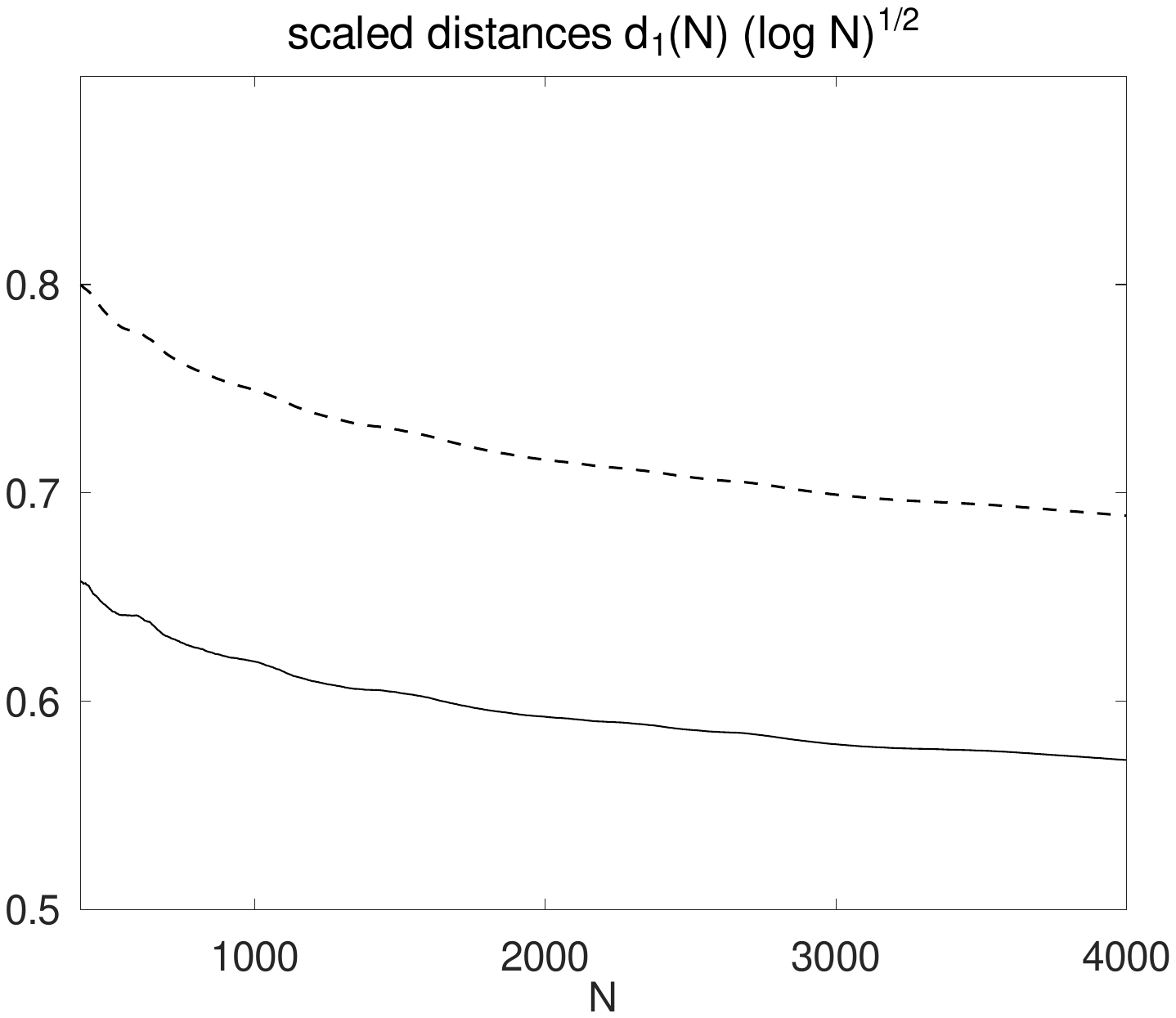}
\end{minipage} 
\hspace{8mm} 
\begin{minipage}{.44\linewidth} 
\includegraphics[scale=.28,angle=0,clip]{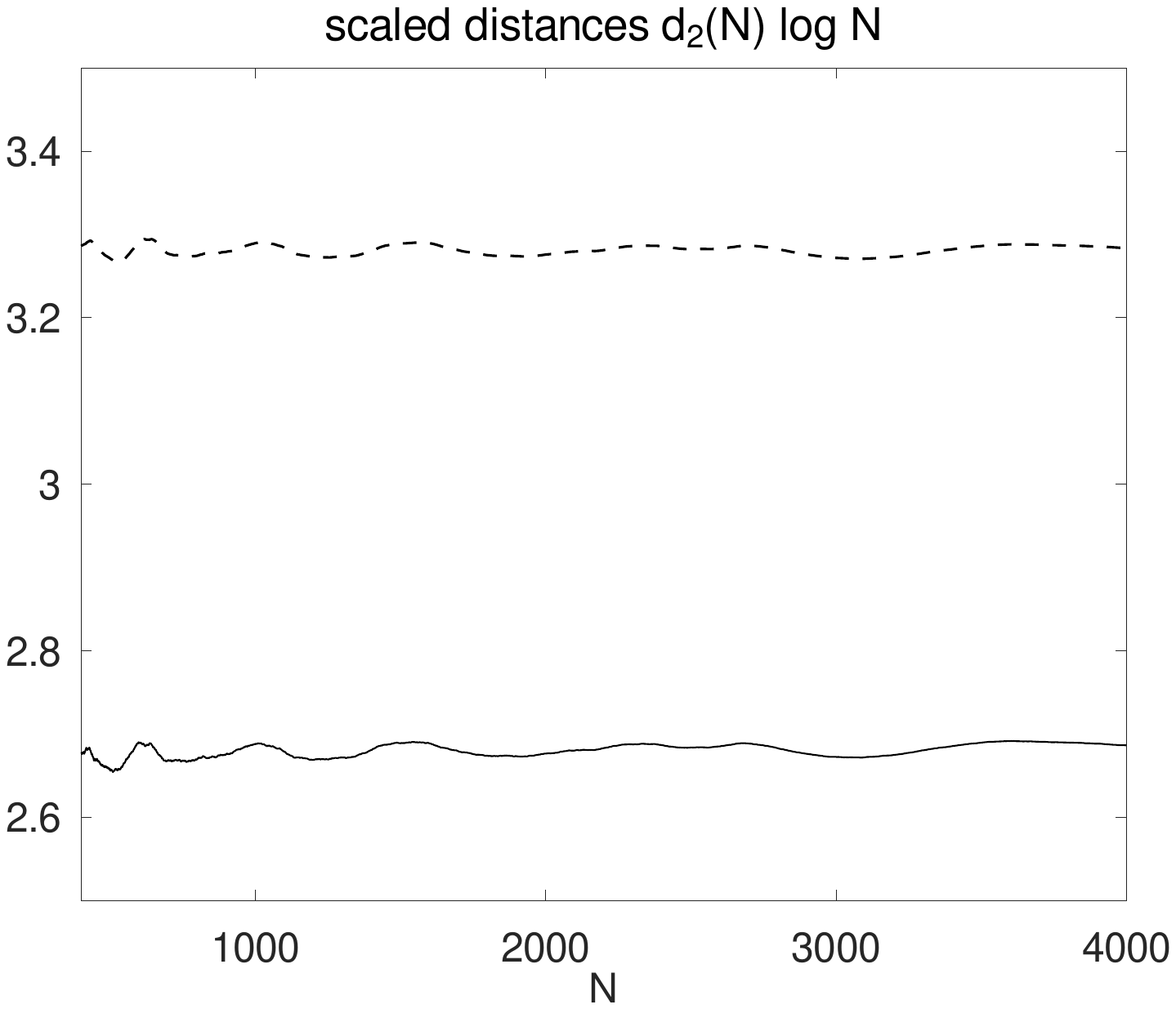}
\end{minipage} \hfill
\newline

\vspace*{-2mm}
Figure 5. \small{$d_1(N) \sqrt{\log N}$ (left)  and $d_2(N) \log N$ (right)  plotted vs.~$N \in [400,4000]$.  Solid: coefficients $\nu_{n,N}$; dashed: coefficients $\lambda_{n,N}$.}
\normalsize

\end{document}